\documentclass[a4paper,12pt]{article}
\usepackage{amssymb}
\usepackage{amsmath}
\usepackage{amsthm}
\usepackage[mathscr]{eucal}
\newtheorem{thm}{Theorem}[section]
\newtheorem{lem}[thm]{Lemma}
\newtheorem{cor}[thm]{Corollary}
\newtheorem{prop}[thm]{Proposition}
\newtheorem{exmp}[thm]{Example}
\makeatletter
\@addtoreset{equation}{section}

\makeatother
\begin{document}
\title{Almost central involutions in split extensions of Coxeter groups by graph automorphisms}
\date{}
\author{Koji Nuida}
\maketitle
\begin{abstract}

In this paper, given a split extension of an arbitrary Coxeter group by automorphisms of the Coxeter graph, we determine the involutions in that extension whose centralizer has finite index.
Our result has applications to many problems such as the isomorphism problem of general Coxeter groups.
In the argument, some properties of certain special elements and of the fixed-point subgroups by graph automorphisms in Coxeter groups, which are of independent interest, are also given.

\end{abstract}

%%%%%%%%%%%%%%%%%%%%%%%%%%%%%%%%%
\section{Introduction}

Let $(W,S)$ be an arbitrary Coxeter system, possibly of infinite rank, and $G$ a group acting on $W$.
We assume that the action of $G$ preserves the set $S$; namely, each element of $G$ gives rise to an automorphism of the Coxeter graph of $(W,S)$.
The subject of this paper is the almost central involutions in the semidirect product $W \rtimes G$ corresponding the action of $G$; that is, involutions which is central in some subgroup of $W \rtimes G$ of finite index.
We determine those involutions in $W \rtimes G$, hence the subgroup generated by those involutions, in terms of the structure of the Coxeter system $(W,S)$ and the action of $G$ on $W$ (Theorem \ref{thm:maintheorem}).
Actually, this subgroup is the product of some finite irreducible components of $W$, specified in terms of the action of $G$, and a subgroup of $G$.
Note that this subgroup is determined by the group structure of $W \rtimes G$ only, so our result can extract some information on the Coxeter group $W$ from the group structure of $W \rtimes G$.
Moreover, if $W \rtimes G$ admits another expression $W' \rtimes G'$ of this type, our result exhibits some relation between the Coxeter groups $W$ and $W'$ through the subgroup in problem (Theorem \ref{thm:cor_maintheorem}).

The main motive of this research is an application to the isomorphism problem of general Coxeter groups; that is, the problem of deciding which Coxeter groups are isomorphic as abstract groups.
An important phase of the problem is to determine whether a given group isomorphism $f$ between two Coxeter groups $W$ and $W'$ maps the reflections in $W$ onto those of $W'$.
As summarized in Section \ref{sec:maintheorem_application}, it is shown by a result of the author's preceding paper \cite{Nui_centra} that both the centralizer of a reflection $t$ in $W$ and that of $f(t)$ in $W'$ are semidirect products satisfying the hypothesis of our main theorem.
Since those centralizers are isomorphic via $f$, our main theorem can derive some properties of $f(t)$ from those of $W$ and of $t$.
In particular, $f(t)$ is a reflection in $W'$ whenever $W$ and $t$ satisfy a certain condition which is independent on the choice of $W'$ and $f$ (Theorem \ref{thm:conditiontobereflection}).
When the condition is actually satisfied will be investigated in a forthcoming paper \cite{Nui_refindep} of the author.
Note that this argument works without any assumption on finiteness of ranks of $W$ or of $W'$, in contrast with most of the preceding results on the isomorphism problem which covers the case of finite ranks only.

For other applications, our result implies that the product of all finite irreducible components of a Coxeter group $W$ is independent on the choice of the generating set $S$ of $W$ (Example \ref{ex:finitepart}).
On the other hand, regarding certain semidirect product decompositions of $W$ into two Coxeter groups which arise from the partition of $S$ into conjugacy classes, our result shows that, under a certain condition, the normal factor possesses no finite irreducible component (Example \ref{ex:semidirectproductdecomposition}).
See Section \ref{sec:maintheorem_example} for further examples.

This paper is organized as follows.
Section \ref{sec:preliminaries} is a preliminary for basics and further remarks on abstract groups and Coxeter groups.
Section \ref{sec:maintheoremandapplications} summarizes the main result and its applications mentioned above.
In Section \ref{sec:essentialelements}, we recall the notion of essential elements in Coxeter groups introduced by Daan Krammer \cite{Kra}, and summarize some properties studied by Krammer and by Luis Paris \cite{Par}.
In Section \ref{sec:fixedpoint}, we give some results on the fixed-point subgroup of a Coxeter group by an automorphism of the Coxeter graph, together with preceding results given by Robert Steinberg \cite{Ste}, by Bernhard M\"{u}hlherr \cite{Muh} and by Masayuki Nanba \cite{Nan}.
Finally, Section \ref{sec:proof_maintheorem} is devoted to the proof of the main theorem.

\noindent
\textbf{Acknowledgement.}
The author would like to express his deep gratitude to everyone who helped him, especially to his supervisor Itaru Terada and also to Kazuhiko Koike for their precious advice and encouragement.
The author had been supported by JSPS Research Fellowship throughout this research.
%%%%
\section{Preliminaries}
\label{sec:preliminaries}
%%%%%
\subsection{On abstract groups}
\label{sec:abstractgroups}
In this subsection, we fix notations for abstract groups, and give some definitions and facts.
Let $G$ be an arbitrary group.
We denote $H \leq G$ if $H$ is a subgroup of $G$, and $H \unlhd G$ if $H$ is a normal subgroup of $G$.
For a subset $X \subseteq G$, let $\langle X \rangle$ and $\langle X \rangle_{\lhd G}$ denote the subgroup and the normal subgroup, respectively, of $G$ generated by $X$.
Put
\[
Z_H(X)=\{g \in H \mid gx=xg \textrm{ for all } x \in X\} \textrm{ for } H \leq G,
\]
so $Z_G(X)$ is the centralizer of $X$ in $G$.
Write
\[
x^g=g^{-1}xg \textrm{ and } X^g=\{x^g \mid x \in X\} \textrm{ for } g,x \in G \textrm{ and } X \subseteq G.
\]
For $H \leq G$, put
\[
\mathrm{Core}_GH=\bigcap_{g \in G}H^g,
\]
the \emph{core} of $H$ in $G$.
It is easily verified that $\mathrm{Core}_GH$ is the unique largest normal subgroup of $G$ contained in $H$.
\begin{lem}
\label{lem:Zofnormalsubgroup}
Let $G$ be a group.
\begin{enumerate}
\item If $H \unlhd G$, then $Z_G(H) \unlhd G$.
\item If $X \subseteq G$, then $Z_G(\langle X \rangle_{\lhd G})=\mathrm{Core}_GZ_G(X)$.
\end{enumerate}
\end{lem}
\begin{proof}
The proof of (1) is straightforward.
For (2), the inclusion $\subseteq$ follows from (1) since $\langle X \rangle \subseteq \langle X \rangle_{\lhd G}$, so it suffices to show that $H \subseteq Z_G(\langle X \rangle_{\lhd G})$ whenever $H \unlhd G$ and $H \subseteq Z_G(X)$.
Now we have $X \subseteq Z_G(H) \unlhd G$ by (1), so $\langle X \rangle_{\lhd G} \subseteq Z_G(H)$, proving the claim.
\end{proof}
Let $\left[G:H\right]$ denote the index of a subgroup $H \leq G$ in $G$.
Recall the following well-known properties:
\begin{eqnarray}
\label{eq:index_chainrule}
&&\textrm{if } G \geq H_1 \geq H_2, \textrm{ then } \left[G:H_2\right]=\left[G:H_1\right]\left[H_1:H_2\right];\\
\label{eq:index_intersection}
&&\textrm{if } H_1,H_2 \leq G, \textrm{ then } \left[G:H_1\right] \geq \left[H_2:H_1 \cap H_2\right].
\end{eqnarray}
From these properties it is easy to deduce that
\begin{eqnarray}
\nonumber
&&\textrm{if } H_1,H_2 \leq G \textrm{ and } \left[G:H_2\right]<\infty, \textrm{ then the followings are equivalent:}\\
\label{eq:index_diamond}
&&\left[G:H_1\right]<\infty;\quad \left[G:H_1 \cap H_2\right]<\infty;\quad \left[H_2:H_1 \cap H_2\right]<\infty.
\end{eqnarray}
\begin{lem}
\label{lem:indexofcore}
Let $H \leq G$.
Then $\left[G:H\right]<\infty$ if and only if $\left[G:\mathrm{Core}_GH\right]<\infty$.
\end{lem}
\begin{proof}
The only nontrivial part is the ``only if'' part.
Let $G=\bigsqcup_{i=1}^nHg_i$ (where $n=\left[G:H\right]<\infty$) be a decomposition into cosets.
Then $\mathrm{Core}_GH=\bigcap_{i=1}^nH^{g_i}$.
Now for $1 \leq k \leq n$, two subgroups $H^{g_k}$ and $H$ have the same (finite) index in $G$, so the subgroup $\bigcap_{i=1}^kH^{g_i}$ has finite index in $\bigcap_{i=1}^{k-1}H^{g_i}$ by (\ref{eq:index_intersection}).
Now iterative use of (\ref{eq:index_chainrule}) yields the desired conclusion.
\end{proof}
We say that an element $g \in G$ is \emph{almost central} in $G$ if $\left[G:Z_G(g)\right]<\infty$.
\begin{cor}
\label{cor:indexandZ}
Let $G$ be a group and $g \in G$.
\begin{enumerate}
\item We have $\left[G:Z_G(\langle g \rangle_{\lhd G})\right]<\infty$ if and only if $g$ is almost central in $G$.
\item If $g$ is almost central in $G$, then all $h \in \langle g \rangle_{\lhd G}$ are almost central in $G$.
\end{enumerate}
\end{cor}
\begin{proof}
The claim (1) follows immediately from Lemmas \ref{lem:Zofnormalsubgroup} (2) and \ref{lem:indexofcore}, and (2) is a consequence of (1) and the observation $Z_G(h) \geq Z_G(\langle g \rangle_{\lhd G})$.
\end{proof}
\begin{lem}
\label{lem:indexinsemidirect}
Let $G_1 \rtimes G_2$ be a semidirect product of two groups, and suppose that $H_i \leq G_i$ has finite index in $G_i$ for $i=1,2$.
Then $\left[G_1 \rtimes G_2:H_1H_2\right]<\infty$.
\end{lem}
\begin{proof}
Decompose $G_i$ as $\bigsqcup_{j=1}^{r_i}g_{i,j}H_i$, where $r_i<\infty$.
Then
\[
G_1 \rtimes G_2=\bigcup_{1 \leq j \leq r_1,\,1 \leq k \leq r_2}g_{1,j}H_1g_{2,k}H_2=\bigcup_{j,k}g_{1,j}g_{2,k}H_1^{g_{2,k}}H_2.
\]
Since $\left[G_1:H_1\right]<\infty$, we have $\left[H_1^{g_{2,k}}:H_1^{g_{2,k}} \cap H_1\right]<\infty$ by (\ref{eq:index_intersection}).
Let $H_1^{g_{2,k}}=\bigsqcup_{\ell=1}^{n_k}h_{k,\ell}(H_1^{g_{2,k}} \cap H_1)$ (where $n_k<\infty$) be the corresponding coset decomposition.
Then we have
\[
G_1 \rtimes G_2=\bigcup_{j,k}\bigcup_{\ell=1}^{n_k}g_{1,j}g_{2,k}h_{k,\ell}(H_1^{g_{2,k}} \cap H_1)H_2 \subseteq \bigcup_{j,k,\ell}g_{1,j}g_{2,k}h_{k,\ell}H_1H_2.
\]
where the last union is taken over the finite set of the $(j,k,\ell)$, as desired.
\end{proof}
%%
%%%%%
\subsection{Coxeter groups}
\label{sec:Coxetergroups}
This subsection summarizes some basic definitions and facts for Coxeter groups, which are found in the book \cite{Hum} unless otherwise noticed, and give further results and remarks.
Some more preliminaries focusing into the two topics, essential elements and fixed-point subgroups by Coxeter graph automorphisms, will be given in Sections \ref{sec:essentialelements} and \ref{sec:fixedpoint}.
%%%%%
\subsubsection{Definitions}
\label{sec:Coxetergroups_definition}
A pair $(W,S)$ of a group $W$ and its generating set $S$ is called a \emph{Coxeter system} if $W$ admits the following presentation
\[
W=\langle S \mid (st)^{m_{s,t}}=1 \textrm{ for all } s,t \in S \textrm{ such that } m_{s,t}<\infty \rangle,
\]
where the $m_{s,t} \in \{1,2,\dots\} \cup \{\infty\}$ are symmetric in $s,t \in S$, and $m_{s,t}=1$ if and only if $s=t$.
A group $W$ is called a \emph{Coxeter group} if some $S \subseteq W$ makes $(W,S)$ a Coxeter system.
The cardinality $|S|$ of $S$ is called the \emph{rank} of $(W,S)$ or of $W$, which is \emph{not} assumed to be finite unless otherwise noticed.
Now $m_{s,t}$ coincides with the order of $st \in W$, so the system $(W,S)$ determines uniquely (up to isomorphism) the \emph{Coxeter graph} denoted by $\Gamma$, that is a simple unoriented graph with vertex set $S$ in which every two vertices $s,t \in S$ is joined by an edge with label $m_{s,t}$ if and only if $m_{s,t} \geq 3$.
(By convention, the label `$3$' is usually omitted when drawing a picture.)

An automorphism of the Coxeter graph $\Gamma$ is briefly called a \emph{graph automorphism of $(W,S)$} of \emph{of $W$}.
Let $\mathrm{Aut}\,\Gamma$ denote the set of the graph automorphisms of $W$.
Then $m_{\tau(s),\tau(t)}=m_{s,t}$ for $\tau \in \mathrm{Aut}\,\Gamma$ and $s,t \in S$, so this $\tau$ extends uniquely to an automorphism of the group $W$ denoted also by $\tau$.

For $I \subseteq S$, let $W_I$ denote the \emph{standard parabolic subgroup} $\langle I \rangle$ of $W$ generated by $I$.
A subgroup conjugate to some $W_I$ is called a \emph{parabolic subgroup}.
(In some context, the term ``parabolic subgroups'' signifies the subgroups $W_I$ themselves only.)
Now $(W_I,I)$ is also a Coxeter system, of which the Coxeter graph $\Gamma_I$ is the full subgraph of $\Gamma$ with vertex set $I$.
If $I$ is (the vertex set of) a connected component of $\Gamma$, then $W_I$ is called an \emph{irreducible component} of $(W,S)$ (or of $W$, if the set $S$ is obvious from the context).
If $\Gamma$ is connected, then $(W,S)$ and $W$ are called \emph{irreducible}.
Now $W$ is the (restricted) direct product of the irreducible components; however, each irreducible component is \emph{not} necessarily directly indecomposable as an abstract group.

Regarding the standard parabolic subgroups, it is well known that
\begin{equation}
\label{eq:intersectionofW_I}
\textrm{if } I,J \subseteq S, \textrm{ then } W_I \cap W_J=W_{I \cap J}.
\end{equation}
Then, since each $w \in W$ is a product of a finite number of elements of $S$, it follows that $W$ possesses a unique minimal standard parabolic subgroup containing $w$, called the \emph{standard parabolic closure} of $w$ and denoted here by $\mathrm{SP}(w)$.
Now the \emph{support} $\mathrm{supp}(w) \subseteq S$ of $w \in W$ is defined by
\[
W_{\mathrm{supp}(w)}=\mathrm{SP}(w).
\]
On the other hand, we have the following fact for parabolic subgroups:
\begin{prop}
[See e.g.\ {\cite[Corollary 7]{Fra-How_rank3}}]
\label{prop:intersectionofparabolic}
Let $I,J \subseteq S$ and $w \in W$.
Then $W_I \cap (W_J)^w=(W_K)^u$ for some $K \subseteq I$ and $u \in W_I$.
Moreover, we have $K \neq I$ whenever $W_I \neq (W_J)^w$.
\end{prop}
This proposition denies the existence of an infinite, properly descending sequence $(W_{I_1})^{w_1} \supset (W_{I_2})^{w_2} \supset \cdots$ of parabolic subgroups with $I_1$ finite, since it enables us to choose the $I_i$ inductively as descending properly.
Thus $W$ also possesses a unique minimal parabolic subgroup containing a given $w \in W$, called the \emph{parabolic closure} of $w$ and denoted here by $\mathrm{P}(w)$.

Let $\ell$ denote the length function of $(W,S)$, namely $\ell(w)$ (where $w \in W$) is the minimal length $n$ of an expression $w=s_1 \cdots s_n$ with $s_i \in S$ (so $\ell(w^{-1})=\ell(w)$).
Such an expression of $w$ with $n=\ell(w)$ is called a \emph{reduced expression}.
The following three well-known properties will be used in the arguments below, without references:
\begin{quote}
if $w \in W$ and $s \in S$, then $\ell(ws)=\ell(w) \pm 1$;\\
for $I \subseteq S$, the length function $\ell_I$ of $(W_I,I)$ agrees with $\ell$ on $W_I$;\\
$\mathrm{supp}(w)=\{s_1,\dots,s_n\}$ for any reduced expression $w=s_1 \cdots s_n$.
\end{quote}
\begin{thm}
[Exchange Condition]
\label{thm:exchangecondition}
Let $w=s_1 \cdots s_n \in W$, $s_i \in S$ and $t \in S$ with $\ell(wt)<\ell(w)$.
Then there exists an index $i$ such that $wt=s_1 \cdots \widehat{s_i} \cdots s_n$ ($s_i$ omitted).
\end{thm}
%%
%%%%%
\subsubsection{Geometric representation and root systems}
\label{sec:rootsystem}
Let $V$ denote the geometric representation space of $W$, that is an $\mathbb{R}$-vector space equipped with the basis $\Pi=\{\alpha_s \mid s \in S\}$ and the symmetric bilinear form $\langle \,,\, \rangle$ determined by
\[
\langle \alpha_s, \alpha_t \rangle=-\cos\frac{\pi}{m_{s,t}} \textrm{ if } m_{s,t}<\infty \textrm{ and } \langle \alpha_s, \alpha_t \rangle=-1 \textrm{ if } m_{s,t}=\infty.
\]
$W$ acts faithfully on $V$ by $s \cdot v=v-2\langle \alpha_s, v\rangle \alpha_s$ for $s \in S$ and $v \in V$, making $\langle \,,\, \rangle$ $W$-invariant.
Let $\Phi=W \cdot \Pi$, $\Phi^+=\Phi \cap \mathbb{R}_{\geq 0}\Pi$ and $\Phi^-=-\Phi^+$ denote, respectively, the \emph{root system}, the set of \emph{positive roots} and the set of \emph{negative roots}.
We have $\Phi=\Phi^+ \sqcup \Phi^-$, and $\Phi$ consists of unit vectors with respect to $\langle \,,\, \rangle$.
For any subset $\Psi \subseteq \Phi$ and $w \in W$, write
\[
\Psi^+=\Psi \cap \Phi^+,\ \Psi^-=\Psi \cap \Phi^- \textrm{ and } \Psi\left[w\right]=\{\gamma \in \Psi^+ \mid w \cdot \gamma \in \Phi^-\}.
\]
Then $\ell(w)$ coincides with the cardinality $|\Phi\left[w\right]|$ of $\Phi\left[w\right]$, so $w=1$ if and only if $\Phi\left[w\right]=\emptyset$.
This implies a further property that
\begin{equation}
\label{eq:rootsdeterminew}
\textrm{ for } w,u \in W, \textrm{ we have } w=u \textrm{ if and only if } \Phi\left[w\right]=\Phi\left[u\right].
\end{equation}
For any $v=\sum_{s \in S}c_s\alpha_s \in V$, the \emph{support} $\mathrm{supp}(v) \subseteq S$ of $v$ is defined by
\[
\mathrm{supp}(v)=\{s \in S \mid c_s \neq 0\}.
\]
For $I \subseteq S$, put
\[
\Pi_I=\{\alpha_s \mid s \in I\} \subset \Pi,\ V_I=\mathrm{span}_{\mathbb{R}}\Pi_I \subset V \textrm{ and } \Phi_I=\Phi \cap V_I. 
\]
Then it is well known (see e.g.\ {\cite[Lemma 4]{Fra-How_nearly}}) that
\begin{equation}
\label{eq:rootsystemofparabolic}
\Phi_I=W_I \cdot \Pi_I,
\end{equation}
the root system of a Coxeter system $(W_I,I)$.
Note that $\Phi\left[w\right] \subseteq \Phi_{\mathrm{supp}(w)}$ for $w \in W$.
Moreover, it is well known that for $I \subseteq S$, any $w \in W$ admits a unique decomposition $w=w^Iw_I$ with $w_I \in W_I$ and $\Phi_I\left[\!\right.w^I\left.\!\right]=\emptyset$.
Note that $\Phi\left[w_I\right]=\Phi_I\left[w\right]$.
This implies that
\begin{equation}
\label{eq:supportandrootofw}
\textrm{if } w \in W \textrm{ and } s \in \mathrm{supp}(w), \textrm{ then } s \in \mathrm{supp}(\gamma) \textrm{ for some } \gamma \in \Phi\left[w\right]
\end{equation}
(if this fails, then $\Phi\left[w\right]=\Phi_I\left[w\right]=\Phi\left[w_I\right]$ where $I=\mathrm{supp}(w) \smallsetminus \{s\}$, so $w=w_I \in W_I$ by (\ref{eq:rootsdeterminew}), contradicting the definition of $\mathrm{supp}(w)$).
Now we prepare a technical lemma which will be used in later sections.
\begin{lem}
\label{lem:lemmaforsupport}
Let $1 \neq w \in W$ and $I=\mathrm{supp}(w) \subseteq S$.
\begin{enumerate}
\item Let $\gamma \in \Phi^+$, $J=\mathrm{supp}(\gamma)$ and suppose that $I \cap J=\emptyset$ and $J$ is adjacent to $I$ in the Coxeter graph $\Gamma$.
Then $w \cdot \gamma \in \Phi_{I \cup J}^+ \smallsetminus \Phi_J$.
\item Suppose that $s \in S \smallsetminus I$ is adjacent to $I$ in $\Gamma$.
Then $\mathrm{supp}(w^s)=I \cup \{s\}$.
\item For $i=1,2$, let $1 \neq u_i \in W$, $J_i=\mathrm{supp}(u_i)$ and suppose that $J_i \cap I=\emptyset$ and $J_2$ is adjacent to $I$ in $\Gamma$.
Then $u_1wu_2 \neq w$.
\end{enumerate}
\end{lem}
\begin{proof}
\textbf{(1)} Since the action of $w \in W_I$ leaves the coefficient in $\gamma$ of any $\alpha_s \in \Pi_J$ unchanged, it suffices to show that $w \cdot \gamma \neq \gamma$.
Take $s \in I$ adjacent to $J$, and $\beta \in \Phi_I^+$ such that $s \in \mathrm{supp}(\beta)$ and $w \cdot \beta \in \Phi_I^-$ (see (\ref{eq:supportandrootofw})).
This choice yields that $\langle \beta, \gamma \rangle<0$ and $\langle w \cdot \beta, \gamma \rangle \geq 0$ since $I \cap J=\emptyset$, showing that $w \cdot \gamma \neq \gamma$ since $\langle \,,\, \rangle$ is $W$-invariant.\\
\textbf{(2)} Put $J=\mathrm{supp}(w^s)$.
Then we have $w=(w^s)^s \in W_{J \cup \{s\}}$ and so $I \subseteq J \cup \{s\}$, therefore $I \subseteq J$ since $s \not\in I$.
On the other hand, $ws \cdot \alpha_s=-w \cdot \alpha_s \in \Phi^- \smallsetminus \{-\alpha_s\}$ by (1), so we have $w^s \cdot \alpha_s \in \Phi^-$ and $s \in J$.
Thus we have $I \cup \{s\} \subseteq J$, while $w^s \in W_{I \cup \{s\}}$, proving the claim.\\
\textbf{(3)} Take $s \in J_2$ adjacent to $I$, and $\gamma \in \Phi\left[\!\right.{u_2}^{-1}\left.\!\right] \subseteq \Phi_{J_2}^+$ with $s \in \mathrm{supp}(\gamma)$ (see (\ref{eq:supportandrootofw})), so $\beta={u_2}^{-1} \cdot \gamma$ lies in $\Phi_{J_2}^-$.
Then $w \cdot \beta \in \Phi^-$ since $I \cap J_2=\emptyset$, while $wu_2 \cdot \beta=w \cdot \gamma \in \Phi^+ \smallsetminus \Phi_{S \smallsetminus I}$ by (1) and so $u_1wu_2 \cdot \beta \in \Phi^+$ since $J_1 \cap I=\emptyset$.
Thus we have $u_1wu_2 \neq w$ as desired.
\end{proof}
For $\gamma=w \cdot \alpha_s \in \Phi$, let $s_\gamma=wsw^{-1}$ denote the \emph{reflection} along the root $\gamma$ acting on $V$ by $s_\gamma \cdot v=v-2 \langle \gamma, v \rangle \gamma$ for $v \in V$.
Let
\[
S^W=\bigcup_{w \in W}wSw^{-1}
\]
denote the set of the reflections in $W$, which depends on the set $S$ in general.
\begin{lem}
\label{lem:orbitofrootisinfinite}
Let $W$ be an infinite irreducible Coxeter group.
Then the orbit $W \cdot \gamma \subseteq \Phi$ of any root $\gamma \in \Phi$ is an infinite set.
\end{lem}
The proof of this lemma requires the following two results:
\begin{prop}
[{\cite[proof of Proposition 4.2]{Deo_root}}]
\label{prop:infinitelymanyroot}
Let $W$ be an infinite irreducible Coxeter group of finite rank, and $I \subset S$ a proper subset.
Then $|\Phi \smallsetminus \Phi_I|=\infty$.
\end{prop}
\begin{prop}
[{\cite[Lemma 2.9]{Nui_indec}}]
\label{prop:lemmaforresidue}
Let $w \in W$ and suppose that $I,J \subseteq S$ are disjoint subsets such that $w \cdot \Pi_I=\Pi_I$ and $w \cdot \Pi_J \subseteq \Phi^-$.
Then we have $\Phi_{I \cup J}\left[w\right]=\Phi_{I \cup J}^+ \smallsetminus \Phi_I$.
\end{prop}
\noindent
\textbf{Proof of Lemma \ref{lem:orbitofrootisinfinite}.}
First we show that, for any $\beta \in \Phi^+$, we have $\langle \beta,\alpha_s \rangle<0$ for some $s \in S$.
This is obvious if $|S|=\infty$ (choose $s \in S \smallsetminus \mathrm{supp}(\beta)$ adjacent in the infinite connected graph $\Gamma$ to the finite set $\mathrm{supp}(\beta)$), so suppose that $|S|<\infty$.
Assume contrary that $\langle \beta, \alpha_s \rangle \geq 0$ for all $s \in S$.
Put $I=\{s \in S \mid \langle \beta,\alpha_s \rangle=0\} \neq S$ (note that $\langle \beta,\beta \rangle=1$), so $s_\beta$ fixes $\Phi_I$ pointwise.
Then for any $s \in S \smallsetminus I$, we have $\langle \beta, \alpha_s \rangle>0$ and $s_\beta \cdot \alpha_s=\alpha_s-2\langle \beta, \alpha_s \rangle \beta \in \Phi^-$.
Thus Proposition \ref{prop:lemmaforresidue} implies that $\Phi\left[s_\beta\right]=\Phi^+ \smallsetminus \Phi_I$, which has cardinality $\ell(s_\beta)<\infty$, contradicting Proposition \ref{prop:infinitelymanyroot}.
Hence the claim of this paragraph holds.

For the lemma, we may assume that $\gamma \in \Phi^+$.
Then by taking $s \in S$ with $\langle \gamma, \alpha_s \rangle<0$ and putting $\gamma_1=s \cdot \gamma$, we have $\gamma_1 \neq \gamma$ and $\gamma_1-\gamma \in \mathbb{R}_{\geq 0}\Pi$.
Iterating, we obtain an infinite sequence $\gamma_0=\gamma$, $\gamma_1$, $\gamma_2,\dots$ of distinct positive roots in $W \cdot \gamma$ inductively, proving the claim.\qed\medskip

We also prepare a technical lemma.
\begin{lem}
\label{lem:reflectionincoset}
Let $\beta,\gamma \in \Phi^+$, $I \subseteq S$ and suppose that $\mathrm{supp}(\gamma) \not\subseteq \mathrm{supp}(\beta)$ and $\mathrm{supp}(\gamma) \not\subseteq I$.
Then $s_\gamma \not\in s_\beta W_I$.
\end{lem}
\begin{proof}
Assume contrary that $s_\gamma=s_\beta w$ for some $w \in W_I$.
Then we have $w \cdot \gamma=s_\beta s_\gamma \cdot \gamma=-s_\beta \cdot \gamma$, while $w \cdot \gamma \in \Phi^+$ and $s_\beta \cdot \gamma \in \Phi^+$ by the hypothesis.
This is a contradiction.
\end{proof}
%%
%%%%%
\subsubsection{Finite, affine and hyperbolic Coxeter groups}
\label{sec:hyperbolicCoxetergroups}
The finite irreducible Coxeter groups are completely classified, as summarized in {\cite[Chapter 2]{Hum}}.
If $I \subseteq S$ and $W_I$ is finite, let $w_0(I)$ denote the unique \emph{longest element} of $W_I$, which is an involution and maps $\Pi_I$ onto $-\Pi_I$.
If $W_I$ is irreducible (but not necessarily finite) and $1 \neq w \in W_I$, then we have $I^w=I$ if and only if $W_I$ is finite and $w=w_0(I)$.
This implies the well-known fact that the center $Z(W_I)$ of an arbitrary $W_I$ is an elementary abelian $2$-group.
Moreover, if $W_I$ is finite but not irreducible, then $w_0(I)=w_0(I_1) \dotsm w_0(I_k)$ where $W_{I_1},\dots,W_{I_k}$ are the irreducible components of $W_I$.
It is well known that, if $w \in W_I$ and $\ell(ws)<\ell(w)$ for all $s \in I$, then $W_I$ is finite and $w=w_0(I)$.
\begin{thm}
[{\cite[Theorem A]{Ric}}]
\label{thm:Ric_involution}
For any involution $w \in W$, there is a finite $W_I$ (where $I \subseteq S$) such that $w$ is conjugate to $w_0(I)$ and $w_0(I) \in Z(W_I)$.
\end{thm}
The cases where $|W_I|<\infty$ and $w_0(I) \in Z(W_I)$ are determined as well.

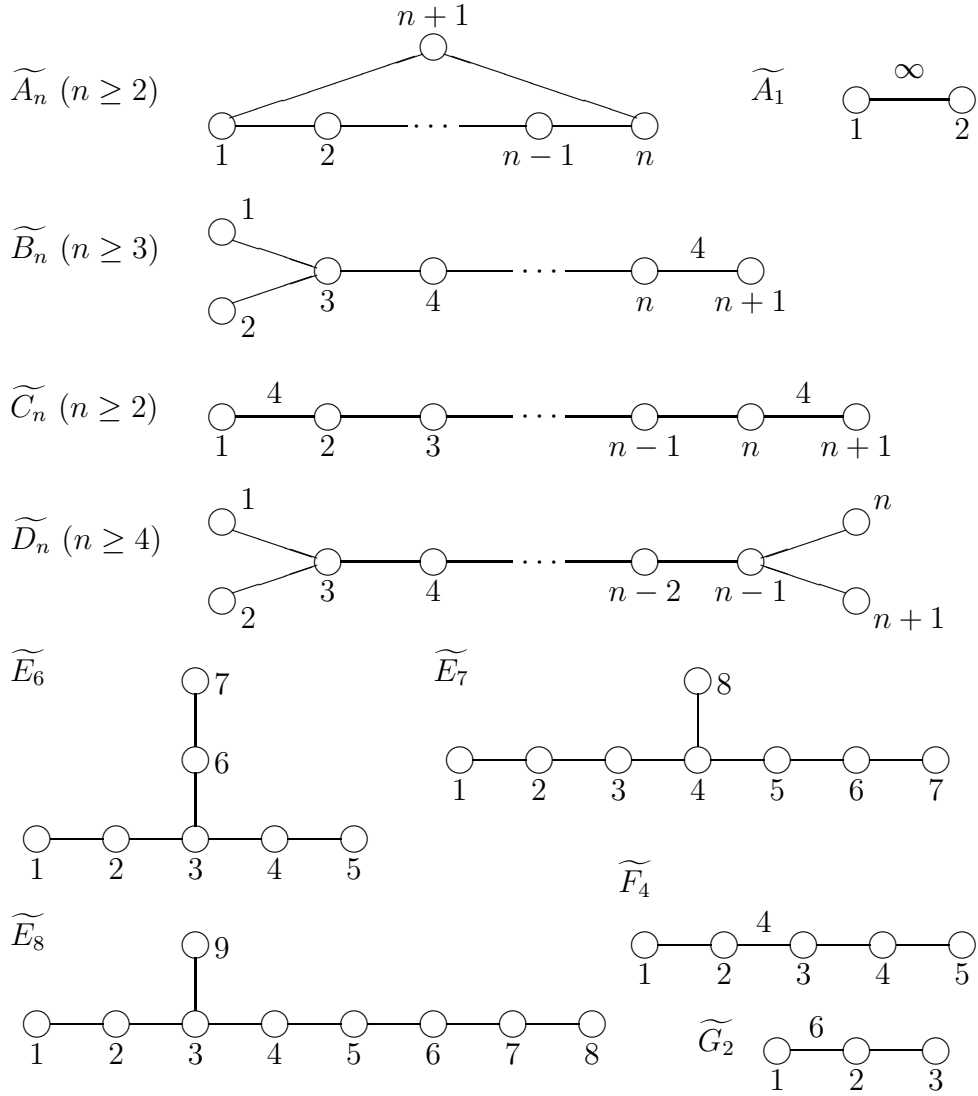
\begin{figure}
\centering
\begin{picture}(370,420)(0,-420)
\put(0,-40){$\widetilde{A_n}$ ($n \geq 2$)}
\put(80,-50){\circle{10}}\put(120,-50){\circle{10}}\put(200,-50){\circle{10}}\put(240,-50){\circle{10}}\put(160,-20){\circle{10}}
\put(83,-47){\line(3,1){73}}\put(237,-47){\line(-3,1){73}}\put(85,-50){\line(1,0){30}}\put(125,-50){\line(1,0){25}}\put(152,-53){$\cdots$}\put(195,-50){\line(-1,0){25}}\put(235,-50){\line(-1,0){30}}
\put(80,-65){\hbox to0pt{\hss$1$\hss}}\put(120,-65){\hbox to0pt{\hss$2$\hss}}\put(200,-65){\hbox to0pt{\hss$n-1$\hss}}\put(240,-65){\hbox to0pt{\hss$n$\hss}}\put(160,-12){\hbox to0pt{\hss$n+1$\hss}}
\put(280,-40){$\widetilde{A_1}$}
\put(320,-40){\circle{10}}\put(360,-40){\circle{10}}
\put(325,-40){\line(1,0){30}}\put(340,-32){\hbox to0pt{\hss$\infty$\hss}}
\put(320,-55){\hbox to0pt{\hss$1$\hss}}\put(360,-55){\hbox to0pt{\hss$2$\hss}}
\put(0,-100){$\widetilde{B_n}$ ($n \geq 3$)}
\put(80,-90){\circle{10}}\put(80,-120){\circle{10}}\put(120,-105){\circle{10}}\put(160,-105){\circle{10}}\put(240,-105){\circle{10}}\put(280,-105){\circle{10}}
\put(84,-93){\line(3,-1){32}}\put(84,-117){\line(3,1){32}}\put(125,-105){\line(1,0){30}}\put(165,-105){\line(1,0){25}}\put(193,-108){$\cdots$}\put(235,-105){\line(-1,0){25}}\put(275,-105){\line(-1,0){30}}\put(260,-100){\hbox to0pt{\hss$4$\hss}}
\put(90,-85){\hbox to0pt{\hss$1$\hss}}\put(90,-130){\hbox to0pt{\hss$2$\hss}}\put(120,-120){\hbox to0pt{\hss$3$\hss}}\put(160,-120){\hbox to0pt{\hss$4$\hss}}\put(240,-120){\hbox to0pt{\hss$n$\hss}}\put(280,-120){\hbox to0pt{\hss$n+1$\hss}}
\put(0,-160){$\widetilde{C_n}$ ($n \geq 2$)}
\put(80,-160){\circle{10}}\put(120,-160){\circle{10}}\put(160,-160){\circle{10}}\put(240,-160){\circle{10}}\put(280,-160){\circle{10}}\put(320,-160){\circle{10}}
\put(85,-160){\line(1,0){30}}\put(125,-160){\line(1,0){30}}\put(165,-160){\line(1,0){25}}\put(193,-163){$\cdots$}\put(235,-160){\line(-1,0){25}}\put(275,-160){\line(-1,0){30}}\put(315,-160){\line(-1,0){30}}\put(100,-155){\hbox to0pt{\hss$4$\hss}}\put(300,-155){\hbox to0pt{\hss$4$\hss}}
\put(80,-175){\hbox to0pt{\hss$1$\hss}}\put(120,-175){\hbox to0pt{\hss$2$\hss}}\put(160,-175){\hbox to0pt{\hss$3$\hss}}\put(240,-175){\hbox to0pt{\hss$n-1$\hss}}\put(280,-175){\hbox to0pt{\hss$n$\hss}}\put(320,-175){\hbox to0pt{\hss$n+1$\hss}}
\put(0,-210){$\widetilde{D_n}$ ($n \geq 4$)}
\put(80,-200){\circle{10}}\put(80,-230){\circle{10}}\put(120,-215){\circle{10}}\put(160,-215){\circle{10}}\put(240,-215){\circle{10}}\put(280,-215){\circle{10}}\put(320,-200){\circle{10}}\put(320,-230){\circle{10}}
\put(84,-203){\line(3,-1){32}}\put(84,-227){\line(3,1){32}}\put(125,-215){\line(1,0){30}}\put(165,-215){\line(1,0){25}}\put(193,-218){$\cdots$}\put(235,-215){\line(-1,0){25}}\put(275,-215){\line(-1,0){30}}\put(316,-203){\line(-3,-1){32}}\put(316,-227){\line(-3,1){32}}
\put(90,-195){\hbox to0pt{\hss$1$\hss}}\put(90,-240){\hbox to0pt{\hss$2$\hss}}\put(120,-230){\hbox to0pt{\hss$3$\hss}}\put(160,-230){\hbox to0pt{\hss$4$\hss}}\put(240,-230){\hbox to0pt{\hss$n-2$\hss}}\put(280,-230){\hbox to0pt{\hss$n-1$\hss}}\put(330,-195){\hbox to0pt{\hss$n$\hss}}\put(340,-240){\hbox to0pt{\hss$n+1$\hss}}
\put(0,-260){$\widetilde{E_6}$}
\put(10,-320){\circle{10}}\put(40,-320){\circle{10}}\put(70,-320){\circle{10}}\put(100,-320){\circle{10}}\put(130,-320){\circle{10}}\put(70,-290){\circle{10}}\put(70,-260){\circle{10}}
\put(15,-320){\line(1,0){20}}\put(45,-320){\line(1,0){20}}\put(75,-320){\line(1,0){20}}\put(105,-320){\line(1,0){20}}\put(70,-315){\line(0,1){20}}\put(70,-285){\line(0,1){20}}
\put(10,-335){\hbox to0pt{\hss$1$\hss}}\put(40,-335){\hbox to0pt{\hss$2$\hss}}\put(70,-335){\hbox to0pt{\hss$3$\hss}}\put(100,-335){\hbox to0pt{\hss$4$\hss}}\put(130,-335){\hbox to0pt{\hss$5$\hss}}\put(80,-295){\hbox to0pt{\hss$6$\hss}}\put(80,-265){\hbox to0pt{\hss$7$\hss}}
\put(160,-260){$\widetilde{E_7}$}
\put(170,-290){\circle{10}}\put(200,-290){\circle{10}}\put(230,-290){\circle{10}}\put(260,-290){\circle{10}}\put(290,-290){\circle{10}}\put(320,-290){\circle{10}}\put(350,-290){\circle{10}}\put(260,-260){\circle{10}}
\put(175,-290){\line(1,0){20}}\put(205,-290){\line(1,0){20}}\put(235,-290){\line(1,0){20}}\put(265,-290){\line(1,0){20}}\put(295,-290){\line(1,0){20}}\put(325,-290){\line(1,0){20}}\put(260,-285){\line(0,1){20}}
\put(170,-305){\hbox to0pt{\hss$1$\hss}}\put(200,-305){\hbox to0pt{\hss$2$\hss}}\put(230,-305){\hbox to0pt{\hss$3$\hss}}\put(260,-305){\hbox to0pt{\hss$4$\hss}}\put(290,-305){\hbox to0pt{\hss$5$\hss}}\put(320,-305){\hbox to0pt{\hss$6$\hss}}\put(350,-305){\hbox to0pt{\hss$7$\hss}}\put(270,-265){\hbox to0pt{\hss$8$\hss}}
\put(0,-360){$\widetilde{E_8}$}
\put(10,-390){\circle{10}}\put(40,-390){\circle{10}}\put(70,-390){\circle{10}}\put(100,-390){\circle{10}}\put(130,-390){\circle{10}}\put(160,-390){\circle{10}}\put(190,-390){\circle{10}}\put(220,-390){\circle{10}}\put(70,-360){\circle{10}}
\put(15,-390){\line(1,0){20}}\put(45,-390){\line(1,0){20}}\put(75,-390){\line(1,0){20}}\put(105,-390){\line(1,0){20}}\put(135,-390){\line(1,0){20}}\put(165,-390){\line(1,0){20}}\put(195,-390){\line(1,0){20}}\put(70,-385){\line(0,1){20}}
\put(10,-405){\hbox to0pt{\hss$1$\hss}}\put(40,-405){\hbox to0pt{\hss$2$\hss}}\put(70,-405){\hbox to0pt{\hss$3$\hss}}\put(100,-405){\hbox to0pt{\hss$4$\hss}}\put(130,-405){\hbox to0pt{\hss$5$\hss}}\put(160,-405){\hbox to0pt{\hss$6$\hss}}\put(190,-405){\hbox to0pt{\hss$7$\hss}}\put(220,-405){\hbox to0pt{\hss$8$\hss}}\put(80,-365){\hbox to0pt{\hss$9$\hss}}
\put(230,-340){$\widetilde{F_4}$}
\put(240,-360){\circle{10}}\put(270,-360){\circle{10}}\put(300,-360){\circle{10}}\put(330,-360){\circle{10}}\put(360,-360){\circle{10}}
\put(245,-360){\line(1,0){20}}\put(275,-360){\line(1,0){20}}\put(305,-360){\line(1,0){20}}\put(335,-360){\line(1,0){20}}\put(285,-355){\hbox to0pt{\hss$4$\hss}}
\put(240,-375){\hbox to0pt{\hss$1$\hss}}\put(270,-375){\hbox to0pt{\hss$2$\hss}}\put(300,-375){\hbox to0pt{\hss$3$\hss}}\put(330,-375){\hbox to0pt{\hss$4$\hss}}\put(360,-375){\hbox to0pt{\hss$5$\hss}}
\put(260,-400){$\widetilde{G_2}$}
\put(290,-400){\circle{10}}\put(320,-400){\circle{10}}\put(350,-400){\circle{10}}
\put(295,-400){\line(1,0){20}}\put(325,-400){\line(1,0){20}}\put(305,-395){\hbox to0pt{\hss$6$\hss}}
\put(290,-415){\hbox to0pt{\hss$1$\hss}}\put(320,-415){\hbox to0pt{\hss$2$\hss}}\put(350,-415){\hbox to0pt{\hss$3$\hss}}
\end{picture}
\caption{List of the affine Coxeter groups}
\label{fig:affine}
\end{figure}

\begin{figure}
\centering
\begin{picture}(360,330)(0,-330)
\put(0,-10){$X_1$}
\put(40,-10){\circle{10}}\put(40,-40){\circle{10}}\put(70,-40){\circle{10}}\put(100,-40){\circle{10}}\put(100,-10){\circle{10}}
\put(45,-10){\line(1,0){50}}\put(40,-15){\line(0,-1){20}}\put(100,-15){\line(0,-1){20}}\put(45,-40){\line(1,0){20}}\put(75,-40){\line(1,0){20}}\put(70,-5){\hbox to0pt{\hss$4$\hss}}
\put(30,-10){\hbox to0pt{\hss$1$\hss}}\put(40,-55){\hbox to0pt{\hss$2$\hss}}\put(70,-55){\hbox to0pt{\hss$3$\hss}}\put(100,-55){\hbox to0pt{\hss$4$\hss}}\put(110,-10){\hbox to0pt{\hss$5$\hss}}
\put(170,-10){$X_2(m_1,m_2)$}\put(150,-40){$\Biggl(\hspace*{90pt}\Biggr)$}\put(160,-30){$3 \leq m_1,m_2 \leq 5$}\put(160,-50){$(m_1,m_2) \neq (3,3)$}
\put(280,-10){\circle{10}}\put(280,-50){\circle{10}}\put(320,-50){\circle{10}}\put(320,-10){\circle{10}}
\put(285,-10){\line(1,0){30}}\put(285,-50){\line(1,0){30}}\put(280,-15){\line(0,-1){30}}\put(320,-15){\line(0,-1){30}}\put(300,-5){\hbox to0pt{\hss$m_1$\hss}}\put(300,-60){\hbox to0pt{\hss$m_2$\hss}}
\put(270,-10){\hbox to0pt{\hss$1$\hss}}\put(270,-50){\hbox to0pt{\hss$2$\hss}}\put(330,-50){\hbox to0pt{\hss$3$\hss}}\put(330,-10){\hbox to0pt{\hss$4$\hss}}
\put(25,-90){$X_3(m_1,m_2,m_3)$}\put(0,-120){$\Biggl(\hspace*{120pt}\Biggr)$}\put(15,-110){$3 \leq m_1,m_2,m_3<\infty$}\put(10,-130){$(m_1,m_2,m_3) \neq (3,3,3)$}
\put(190,-100){\circle{10}}\put(160,-130){\circle{10}}\put(220,-130){\circle{10}}
\put(186,-103){\line(-1,-1){23}}\put(194,-103){\line(1,-1){23}}\put(165,-130){\line(1,0){50}}\put(165,-110){\hbox to0pt{\hss$m_1$\hss}}\put(215,-110){\hbox to0pt{\hss$m_3$\hss}}\put(190,-140){\hbox to0pt{\hss$m_2$\hss}}
\put(190,-90){\hbox to0pt{\hss$1$\hss}}\put(150,-140){\hbox to0pt{\hss$2$\hss}}\put(230,-140){\hbox to0pt{\hss$3$\hss}}
\put(260,-100){$Y_1$}
\put(270,-120){\circle{10}}\put(310,-120){\circle{10}}\put(350,-100){\circle{10}}\put(350,-140){\circle{10}}
\put(275,-120){\line(1,0){30}}\put(314,-117){\line(2,1){31}}\put(314,-123){\line(2,-1){31}}\put(290,-115){\hbox to0pt{\hss$5$\hss}}
\put(270,-135){\hbox to0pt{\hss$1$\hss}}\put(310,-135){\hbox to0pt{\hss$2$\hss}}\put(360,-100){\hbox to0pt{\hss$3$\hss}}\put(360,-145){\hbox to0pt{\hss$4$\hss}}
\put(0,-160){$Y_2$}
\put(00,-190){\circle{10}}\put(40,-190){\circle{10}}\put(80,-190){\circle{10}}\put(120,-170){\circle{10}}\put(120,-210){\circle{10}}
\put(5,-190){\line(1,0){30}}\put(45,-190){\line(1,0){30}}\put(84,-187){\line(2,1){31}}\put(84,-193){\line(2,-1){31}}\put(20,-185){\hbox to0pt{\hss$5$\hss}}
\put(0,-205){\hbox to0pt{\hss$1$\hss}}\put(40,-205){\hbox to0pt{\hss$2$\hss}}\put(80,-205){\hbox to0pt{\hss$3$\hss}}\put(130,-170){\hbox to0pt{\hss$4$\hss}}\put(130,-215){\hbox to0pt{\hss$5$\hss}}
\put(165,-170){$Y_3(m)$}\put(150,-190){($3 \leq m \leq 5$)}
\put(230,-180){\circle{10}}\put(260,-180){\circle{10}}\put(290,-180){\circle{10}}\put(320,-180){\circle{10}}\put(350,-180){\circle{10}}
\put(235,-180){\line(1,0){20}}\put(265,-180){\line(1,0){20}}\put(295,-180){\line(1,0){20}}\put(325,-180){\line(1,0){20}}\put(245,-175){\hbox to0pt{\hss$5$\hss}}\put(335,-175){\hbox to0pt{\hss$m$\hss}}
\put(230,-195){\hbox to0pt{\hss$1$\hss}}\put(260,-195){\hbox to0pt{\hss$2$\hss}}\put(290,-195){\hbox to0pt{\hss$3$\hss}}\put(320,-195){\hbox to0pt{\hss$4$\hss}}\put(350,-195){\hbox to0pt{\hss$5$\hss}}
\put(10,-240){$Y_4(m)$}\put(0,-260){($m=4,5$)}
\put(70,-250){\circle{10}}\put(100,-250){\circle{10}}\put(130,-250){\circle{10}}\put(160,-250){\circle{10}}
\put(75,-250){\line(1,0){20}}\put(105,-250){\line(1,0){20}}\put(135,-250){\line(1,0){20}}\put(85,-245){\hbox to0pt{\hss$5$\hss}}\put(145,-245){\hbox to0pt{\hss$m$\hss}}
\put(70,-265){\hbox to0pt{\hss$1$\hss}}\put(100,-265){\hbox to0pt{\hss$2$\hss}}\put(130,-265){\hbox to0pt{\hss$3$\hss}}\put(160,-265){\hbox to0pt{\hss$4$\hss}}
\put(210,-240){$Y_5$}
\put(240,-250){\circle{10}}\put(270,-250){\circle{10}}\put(300,-250){\circle{10}}\put(330,-250){\circle{10}}
\put(245,-250){\line(1,0){20}}\put(275,-250){\line(1,0){20}}\put(305,-250){\line(1,0){20}}\put(285,-245){\hbox to0pt{\hss$5$\hss}}
\put(240,-265){\hbox to0pt{\hss$1$\hss}}\put(270,-265){\hbox to0pt{\hss$2$\hss}}\put(300,-265){\hbox to0pt{\hss$3$\hss}}\put(330,-265){\hbox to0pt{\hss$4$\hss}}
\put(0,-310){$Y_6(m_1,m_2)$}\put(60,-310){$\Biggl(\hspace*{140pt}\Biggr)$}\put(70,-300){$3 \leq m_2 \leq m_1<\infty$,\ $5 \leq m_1$}\put(80,-320){$(m_1,m_2) \neq (5,3),(6,3)$}
\put(240,-310){\circle{10}}\put(280,-310){\circle{10}}\put(320,-310){\circle{10}}
\put(245,-310){\line(1,0){30}}\put(285,-310){\line(1,0){30}}\put(260,-305){\hbox to0pt{\hss$m_1$\hss}}\put(300,-305){\hbox to0pt{\hss$m_2$\hss}}
\put(240,-325){\hbox to0pt{\hss$1$\hss}}\put(280,-325){\hbox to0pt{\hss$2$\hss}}\put(320,-325){\hbox to0pt{\hss$3$\hss}}
\end{picture}
\caption{List of the compact hyperbolic Coxeter groups}
\label{fig:compacthyperbolic}
\end{figure}
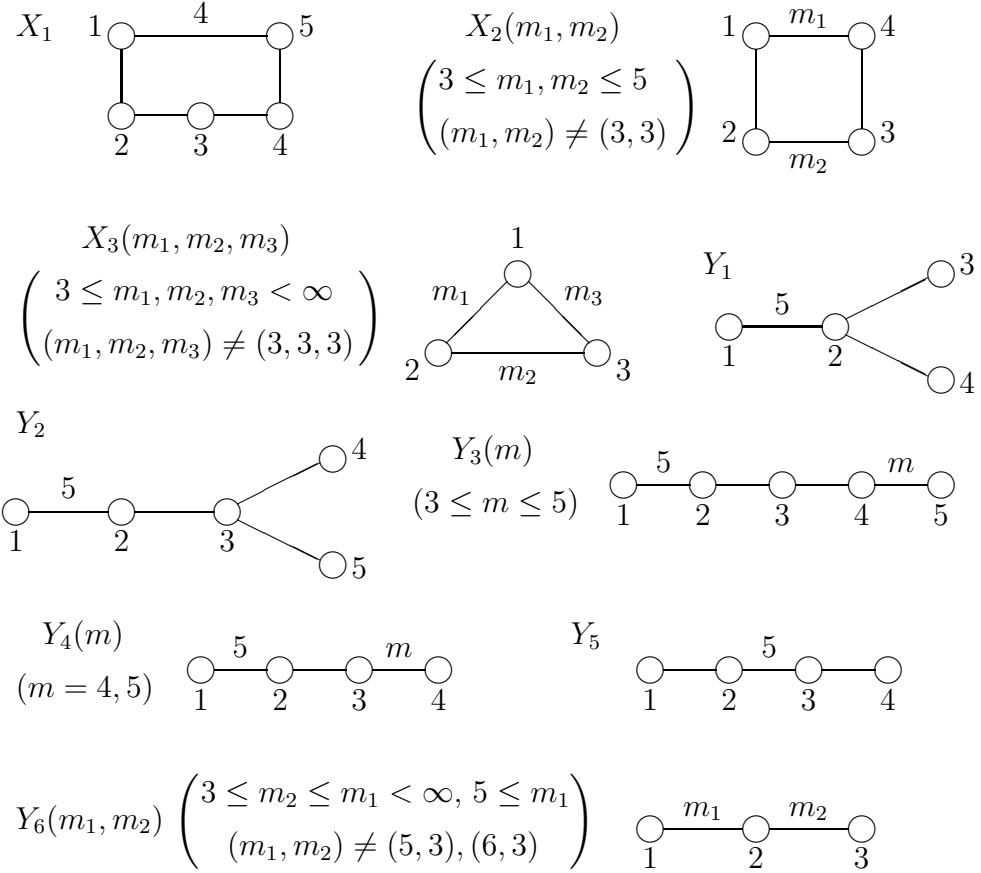

Let $W$ be an irreducible Coxeter group of finite rank.
Then $W$ is called \emph{affine} or \emph{compact hyperbolic}, respectively, if the bilinear form $\langle \,,\, \rangle$ satisfies that (1) it is positive semidefinite or nondegenerate, respectively; (2) it is not positive definite; and (3) its restriction to any proper subspace $V_I \subset V$ (where $I \subset S$) is positive definite.
(See {\cite[Section 6.8]{Hum}} for another definition of compact hyperbolicness and its equivalence to ours.)
The next proposition says that these are the minimal non-finite irreducible Coxeter groups.
\begin{prop}
\label{prop:charofminimalnonfinite}
Let $W$ be a Coxeter group of finite rank.
\begin{enumerate}
\item ({\cite[Theorem 6.4]{Hum}}) We have $|W|<\infty$ if and only if $\langle \,,\, \rangle$ is positive definite.
\item If $|W|=\infty$ and every proper standard parabolic subgroup $W_I \subset W$ is finite, then $W$ is irreducible, and is either affine or compact hyperbolic.
\end{enumerate}
\end{prop}
\begin{proof}
For (2), it is easy to show that this $W$ is irreducible.
Thus by (1) and the definition of compact hyperbolicness, it now suffices to show that this $\langle \,,\, \rangle$ is positive semidefinite if it is degenerate.
This follows from the observation that now $V$ is the sum of a positive definite subspace $V_{S \smallsetminus \{s\}}$ (where $s \in S$; see (1)) of codimension $1$ and the nonzero radical $V^\perp$ of $V$ (note that $V^\perp \not\subseteq V_{S \smallsetminus \{s\}}$).
\end{proof}
The affine and the compact hyperbolic Coxeter groups are completely determined in {\cite[Chapter 2 and Section 6.9]{Hum}}.
See the lists in Figures \ref{fig:affine} and \ref{fig:compacthyperbolic}, where we abbreviate $s_i$ to $i$.
Note that the names of the compact hyperbolic Coxeter groups given here are not standard and are very temporary.
%%%%%
%%%%%
%%%%%

\begin{figure}
\centering
\begin{picture}(390,90)(0,-75)
\put(0,-10){$A_\infty$}
\put(40,-0){\circle{10}}\put(70,-0){\circle{10}}\put(100,-0){\circle{10}}
\put(45,-0){\line(1,0){20}}\put(75,-0){\line(1,0){20}}\put(105,-0){\line(1,0){15}}\put(125,-3){$\cdots$}
\put(40,-15){\hbox to0pt{\hss$1$\hss}}\put(70,-15){\hbox to0pt{\hss$2$\hss}}\put(100,-15){\hbox to0pt{\hss$3$\hss}}
\put(160,-10){$A_{\pm \infty}$}
\put(230,-0){\circle{10}}\put(260,-0){\circle{10}}\put(290,-0){\circle{10}}\put(320,-0){\circle{10}}\put(350,-0){\circle{10}}
\put(195,-3){$\cdots$}\put(225,-0){\line(-1,0){10}}\put(235,-0){\line(1,0){20}}\put(265,-0){\line(1,0){20}}\put(295,-0){\line(1,0){20}}\put(325,-0){\line(1,0){20}}\put(355,-0){\line(1,0){10}}\put(370,-3){$\cdots$}
\put(230,-15){\hbox to0pt{\hss$-2$\hss}}\put(260,-15){\hbox to0pt{\hss$-1$\hss}}\put(290,-15){\hbox to0pt{\hss$0$\hss}}\put(320,-15){\hbox to0pt{\hss$1$\hss}}\put(350,-15){\hbox to0pt{\hss$2$\hss}}
\put(0,-60){$B_\infty$}
\put(40,-50){\circle{10}}\put(70,-50){\circle{10}}\put(100,-50){\circle{10}}\put(130,-50){\circle{10}}
\put(45,-50){\line(1,0){20}}\put(75,-50){\line(1,0){20}}\put(105,-50){\line(1,0){20}}\put(135,-50){\line(1,0){15}}\put(155,-53){$\cdots$}\put(55,-45){\hbox to0pt{\hss$4$\hss}}
\put(40,-65){\hbox to0pt{\hss$1$\hss}}\put(70,-65){\hbox to0pt{\hss$2$\hss}}\put(100,-65){\hbox to0pt{\hss$3$\hss}}\put(130,-65){\hbox to0pt{\hss$4$\hss}}
\put(200,-60){$D_\infty$}
\put(240,-45){\circle{10}}\put(240,-75){\circle{10}}\put(270,-60){\circle{10}}\put(300,-60){\circle{10}}\put(330,-60){\circle{10}}
\put(245,-47){\line(3,-2){19}}\put(245,-73){\line(3,2){19}}\put(275,-60){\line(1,0){20}}\put(305,-60){\line(1,0){20}}\put(335,-60){\line(1,0){15}}\put(355,-63){$\cdots$}
\put(230,-45){\hbox to0pt{\hss$1$\hss}}\put(230,-85){\hbox to0pt{\hss$2$\hss}}\put(270,-75){\hbox to0pt{\hss$3$\hss}}\put(300,-75){\hbox to0pt{\hss$4$\hss}}\put(330,-75){\hbox to0pt{\hss$5$\hss}}
\end{picture}
\caption{Some Coxeter groups of infinite ranks}
\label{fig:nonf.g.Coxetergroups}
\end{figure}
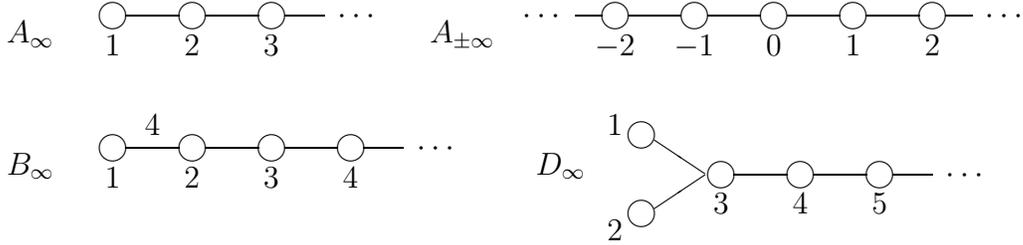

On the other hand, it is shown in {\cite[Proposition 4.14]{Nui_indec}} that the infinite irreducible Coxeter groups of infinite ranks, in which every proper standard parabolic subgroup of finite rank is finite, are exhausted by Figure \ref{fig:nonf.g.Coxetergroups}.
%%%%%
%%%%%
%%%%%
\subsubsection{On centralizers and normalizers in Coxeter groups}
\label{sec:ZinCoxetergroups}
The centralizers and the normalizers in Coxeter groups play important roles in our arguments.
Here we summarize some properties which we require.
\begin{lem}
[e.g.\ {\cite[Lemma 4.4]{Nui_indec}}]
\label{lem:Zoflongestelement}
Let $W_I$ be a finite standard parabolic subgroup of $W$ such that $w_0(I) \in Z(W_I)$.
Then the centralizer $Z_W(w_0(I))$ coincides with the normalizer $N_W(W_I)$ of $W_I$ in $W$.
\end{lem}
\begin{prop}
\label{prop:charoffinitepart}
Let $W$ be an infinite irreducible Coxeter group.
Then no involution in $W$ is almost central in $W$ (see Section \ref{sec:abstractgroups} for terminology).
\end{prop}
\begin{proof}
First, if $s \in S$, then $W$ acts transitively on the conjugacy class of $s$ in $W$, which is an infinite set (Lemma \ref{lem:orbitofrootisinfinite}), so the kernel of this action is $Z_W(s)$ and has infinite index in $W$.
Thus $s$ is not almost central.

By Theorem \ref{thm:Ric_involution}, it suffices to prove that the longest element $w_0(I)$ of any finite $W_I \neq 1$ with $w_0(I) \in Z(W_I)$ is not almost central.
Note that $Z_W(w_0(I))=N_W(W_I)$ (Lemma \ref{lem:Zoflongestelement}), while $\left[N_W(W_I):Z_W(W_I)\right]<\infty$ since $|W_I|<\infty$.
By the first paragraph, $Z_W(s)$ has infinite index in $W$ for any $s \in I$, so do $Z_W(W_I)$ (see (\ref{eq:index_chainrule})) and $Z_W(w_0(I))$, as desired.
\end{proof}
Finally, in {\cite[Theorem 3.1]{Nui_indec}}, the centralizer of a normal subgroup generated by involutions in an irreducible $W$ is completely determined.
The following observation is an easy consequence of the result.
\begin{prop}
[See {\cite[Theorem 3.1]{Nui_indec}}]
\label{prop:Zofnormalinvolutions}
Suppose that $W$ is an arbitrary Coxeter group, and $H \leq W$ is a subgroup generated by involutions which is normal in $W$.
Then $Z_W(H)$ is also generated by involutions.
\end{prop}
%%
%%%%%
\section{The main theorem and its applications}
\label{sec:maintheoremandapplications}
The first subsection of this section summarizes the main theorem of this paper (Theorem \ref{thm:maintheorem}) and its corollary (Theorem \ref{thm:cor_maintheorem}) together with some notational remarks.
The second subsection consists of some examples, and explains what our theorem yields in these cases.
Finally, the third subsection is devoted to an application of our theorem to the analysis of the isomorphism problem of Coxeter groups (the problem of deciding which Coxeter groups are isomorphic as abstract groups), which is the original motive of this research.
%%%
\subsection{Main theorem}
\label{sec:maintheorem}
First we prepare some notations.
Let $W$ be an arbitrary Coxeter group, and $G$ a group acting on $W$ via a map $\rho:G \to \mathrm{Aut}\,\Gamma$, $g \mapsto \rho_g$, yielding the semidirect product $W \rtimes G$ with respect to $\rho$.
Let $\mathcal{C}_W^{\mathrm{fin}}$ and $\mathcal{C}_W^{\mathrm{inf}}$ be the set of the finite and the infinite irreducible components of $W$, respectively, and $\mathcal{C}_W=\mathcal{C}_W^{\mathrm{fin}} \cup \mathcal{C}_W^{\mathrm{inf}}$.
Then the $G$-action permutes the elements of each of $\mathcal{C}_W$, $\mathcal{C}_W^{\mathrm{fin}}$ and $\mathcal{C}_W^{\mathrm{inf}}$.
Let $\rho^\dagger:G \to \mathrm{Sym}(\mathcal{C}_W^{\mathrm{fin}})$, $g \mapsto \rho^\dagger_g$, denote the induced permutation representation of $G$ on $\mathcal{C}_W^{\mathrm{fin}}$.
For $\mathcal{C} \subseteq \mathcal{C}_W$, let $W({\mathcal{C}})$ be the product of the irreducible components in $\mathcal{C}$, and put $W_{\mathrm{fin}}=W(\mathcal{C}_W^{\mathrm{fin}})$ and $W_{\mathrm{inf}}=W(\mathcal{C}_W^{\mathrm{inf}})$.
Moreover, for an arbitrary group $H$, let $H_{\mathrm{ACI}}$ be the set of the almost central involutions in $H$ (see Section \ref{sec:abstractgroups} for the terminology).

Now our main theorem is as follows:
\begin{thm}
\label{thm:maintheorem}
Here we adopt the above notations.
\begin{enumerate}
\item Let $wg$ be an involution in $W \rtimes G$ with $w \in W$ and $g \in G$.
Then $wg$ is almost central in $W \rtimes G$ if and only if $w \in W(\mathcal{O}_\rho)$ and $g \in G_\rho \cup \{1\}$, where $G_\rho$ is the set of all $h \in G_{\mathrm{ACI}}$ satisfying the following condition:
\begin{eqnarray}
\nonumber
&\textrm{$\rho_h$ is identity on all irreducible components of $W$}\\
\label{eq:condition_maintheorem}
&\textrm{except a finite number of finite irreducible components},
\end{eqnarray}
and $\mathcal{O}_{\rho} \subseteq \mathcal{C}_W^{\mathrm{fin}}$ is the union of the $\rho^\dagger(G)$-orbits with finite cardinalities.
\item We have 
\[
\langle (W \rtimes G)_{\mathrm{ACI}} \rangle=W(\mathcal{O}_{\rho}) \rtimes \langle G_\rho \rangle.
\]
\end{enumerate}
\end{thm}
Note that, assuming Theorem \ref{thm:indexoffixedpointsubgroup} below, the condition (\ref{eq:condition_maintheorem}) is equivalent to the finiteness of the index $\left[\!\right.W:W^{\rho_g}\left.\!\right]$ of the fixed-point subgroup $W^{\rho_g}$ by $\rho_g$.
The proof of Theorem \ref{thm:maintheorem} is postponed until Section \ref{sec:proof_maintheorem}.

Since the subgroup $\langle H_{\mathrm{ACI}} \rangle$ of a group $H$ is determined by the isomorphism type of $H$ only, we obtain the following consequence.
\begin{thm}
\label{thm:cor_maintheorem}
For $i=1,2$, let $W_i \rtimes G_i$ be a semidirect product (via $\rho_i:G_i \to \mathrm{Aut}\,\Gamma_i$) as in Theorem \ref{thm:maintheorem}, and $f:W_1 \rtimes G_1 \overset{\sim}{\to} W_2 \rtimes G_2$ a group isomorphism.
Then $f$ maps $W_1(\mathcal{O}_{\rho_1}) \rtimes (G_1)_{\rho_1}$ onto $W_2(\mathcal{O}_{\rho_2}) \rtimes (G_2)_{\rho_2}$.
\end{thm}
%%
%%%%%
\subsection{Examples}
\label{sec:maintheorem_example}
First we observe that, if $|G|<\infty$, then every $\rho^\dagger(G)$-orbit in $\mathcal{C}_W^{\mathrm{fin}}$ is finite, so $\mathcal{O}_\rho=\mathcal{C}_W^{\mathrm{fin}}$ in Theorem \ref{thm:maintheorem}, therefore $\langle (W \rtimes G)_{\mathrm{ACI}} \rangle=W_{\mathrm{fin}} \rtimes G_\rho$ and $G_\rho$ is generated by all involutions $h \in G$ satisfying (\ref{eq:condition_maintheorem}).
\begin{exmp}
\label{ex:finitepart}
Let $W$ be an arbitrary Coxeter group.
Then, by putting $G=1$, Theorem \ref{thm:maintheorem} shows that $\langle W_{\mathrm{ACI}} \rangle=W_{\mathrm{fin}}$.
Thus if $f:W \overset{\sim}{\to} W'$ is a group isomorphism between two Coxeter groups, then $f(W_{\mathrm{fin}})=W'{}_{\mathrm{fin}}$; hence, by taking $W'=W$ and $f=\mathrm{id}_W$, it follows that the factor $W_{\mathrm{fin}}$ is independent on the choice of the generating set $S \subseteq W$.
\end{exmp}
Example \ref{ex:finitepart} is slightly generalized as follows:
\begin{exmp}
\label{ex:wreathproduct}
Let $W \,\mathrm{wr}\, S_n=W^n \rtimes S_n$ denote the wreath product of $W$ with the symmetric group $S_n$ on $n$ letters, so $\sigma \in S_n$ acts on $(w_1,\dots,w_n) \in W^n$ by $\rho_\sigma(w_1,\dots,w_n)=(w_{\sigma^{-1}(1)},\dots,w_{\sigma^{-1}(n)})$.
Then Theorem \ref{thm:maintheorem} implies that
\[
\langle (W \,\mathrm{wr}\, S_n)_{\mathrm{ACI}} \rangle=
\begin{cases}
W \,\mathrm{wr}\, S_n & \textrm{if } |W|<\infty;\\
W_{\mathrm{fin}}{}^n & \textrm{if } |W|=\infty.
\end{cases}
\]
Indeed, if $|W|=\infty$, then $W$ possesses either an infinite irreducible component or infinitely many finite irreducible components, so no non-identity $\sigma \in S_n$ satisfies the condition (\ref{eq:condition_maintheorem}) in any case.
\end{exmp}
We say that an irreducible component $W_I$ of $W$ has \emph{finite multiplicity in $W$} if $W$ possesses only finitely many irreducible components with Coxeter graph isomorphic to $\Gamma_I$.
Note that, even if $|G|=\infty$, the factor $W(\mathcal{O}_\rho)$ in the theorem contains all $W_I \in \mathcal{C}_W^{\mathrm{fin}}$ with finite multiplicities.
\begin{exmp}
\label{ex:GisAutGamma}
Let $G=\mathrm{Aut}\,\Gamma * \mathrm{Aut}\,\Gamma$ be the free product of two copies of $\mathrm{Aut}\,\Gamma$, and $\rho:G \to \mathrm{Aut}\,\Gamma$ the map obtained by forgetting the distinction of the two factors $\mathrm{Aut}\,\Gamma$ of $G$.
Then $\mathcal{O}_\rho$ is the set of all $W_I \in \mathcal{C}_W^{\mathrm{fin}}$ with finite multiplicities, and $G_\rho=1$ since we have $G_{\mathrm{ACI}}=\emptyset$ by properties of free products.
Roughly speaking, Theorem \ref{thm:maintheorem} extracts the finite irreducible components of $W$ with finite multiplicities in this manner.
\end{exmp}
For the final example, we prepare some further facts and notations.
Let $\Gamma^{\mathrm{odd}}$ denote the \emph{odd Coxeter graph} of a Coxeter group $W$, which is the subgraph of $\Gamma$ obtained by removing all the edges with non-odd labels.
It is well known (see {\cite[Exercise 5.3]{Hum}}) that two orbits $W \cdot \alpha_s$ and $W \cdot \alpha_t$ (where $s,t \in S$) intersects nontrivially if and only if $s$ and $t$ lie in the same connected component of $\Gamma^{\mathrm{odd}}$.
Let $S=S_1 \sqcup S_2$ be a partition where both factors are unions of connected components of $\Gamma^{\mathrm{odd}}$, and $\Phi_{\sim S_1}=\bigcup_{s \in S_1}W \cdot \alpha_s \subseteq \Phi$.
Now a general theorem of Vinay V.\ Deodhar \cite{Deo_refsub} or of Matthew Dyer \cite{Dye} shows that the subgroup $W(\Phi_{\sim S_1})$ generated by the reflections $s_\gamma$ along $\gamma \in \Phi_{\sim S_1}$, which is normal in $W$ since $\Phi_{\sim S_1}$ is $W$-invariant, is a Coxeter group.
Moreover, the set $\Phi_{\sim S_1}$ plays the role of a root system of $W(\Phi_{\sim S_1})$; for example, any non-identity $w \in W(\Phi_{\sim S_1})$ sends some $\gamma \in \Phi_{\sim S_1}^+$ to a negative root.

Now we show that $W$ decomposes as $W(\Phi_{\sim S_1}) \rtimes W_{S_2}$.
First, if $1 \neq w \in W(\Phi_{\sim S_1}) \cap W_{S_2}$, then $w \cdot \gamma \in \Phi^-$ for some $\gamma \in \Phi_{\sim S_1}^+$ as mentioned above, and $\gamma \in \Phi_{S_2}$ since $w \in W_{S_2}$.
Now by (\ref{eq:rootsystemofparabolic}), we have $\gamma \in W \cdot \alpha_s \cap W \cdot \alpha_t$ for some $s \in S_1$ and $t \in S_2$, contradicting the choice of the partition $S=S_1 \sqcup S_2$.
Thus we have $W(\Phi_{\sim S_1}) \cap W_{S_2}=1$, while $S \subseteq W(\Phi_{\sim S_1})W_{S_2}$ generates $W$, yielding the desired decomposition.

Moreover, this argument also shows that each $s \in S_2$ preserves the set $\Phi_{\sim S_1}^+$ of positive roots of $W(\Phi_{\sim S_1})$ (since $\alpha_s \not\in \Phi_{\sim S_1}$), so also the set of simple roots of $W(\Phi_{\sim S_1})$, therefore $W_{S_2}$ acts on $W(\Phi_{\sim S_1})$ as graph automorphisms.
Thus Theorem \ref{thm:maintheorem} yields the following observation:
\begin{exmp}
\label{ex:semidirectproductdecomposition}
In the situation, suppose further that $W$ is infinite and irreducible.
Then $\langle W_{\mathrm{ACI}} \rangle=1$ (Example \ref{ex:finitepart}), while $W(\Phi_{\sim S_1})(\mathcal{O}_\rho)$ contains all the finite irreducible components of $W(\Phi_{\sim S_1})$ with finite multiplicities as mentioned above.
Since $1=W(\Phi_{\sim S_1})(\mathcal{O}_\rho) \rtimes (W_{S_2})_\rho$ (Theorem \ref{thm:cor_maintheorem}), it follows that no finite irreducible component of $W(\Phi_{\sim S_1})$ has finite multiplicity in $W(\Phi_{\sim S_1})$.\\
\indent
In addition, if $W_{S_2}$ is finite, then we have $W(\Phi_{\sim S_1})(\mathcal{O}_\rho)=W(\Phi_{\sim S_1})_{\mathrm{fin}}$.
Now it follows that $W(\Phi_{\sim S_1})$ possesses no finite irreducible component.
\end{exmp}
%%
%%%%%
\subsection{Application to the isomorphism problem of Coxeter groups}
\label{sec:maintheorem_application}
An important phase of the isomorphism problem of Coxeter groups is of deciding whether a given group isomorphism $f:W \overset{\sim}{\to} W'$ between two Coxeter groups $W$ and $W'$ (with generating sets $S$ and $S'$, respectively) maps the set $S^W$ of reflections in $W$ onto that ${S'}^{W'}$ in $W'$; or, whether the subset $S^W$ of $W$ is independent on the choice of $S$.
Note that, as is shown in {\cite[Lemma 3.7]{Bra-McC-Muh-Neu}}, we have $f(S^W)={S'}^{W'}$ if and only if $f(S) \subseteq {S'}^{W'}$.
Roughly speaking, our result below measures how $f(s)$ differs from reflections for each $s \in S$, within a certain compass.
In most successful cases, the result is able to show that all $f(s)$ are reflections in $W'$ (see Theorem \ref{thm:conditiontobereflection}).

Note that our results cover the case $|S|=\infty$ as well, in contrast with almost all of the preceding results on the isomorphism problem which cover the case of finite ranks only.
%%%%%
\subsubsection{Preliminaries on centralizers and normalizers}
\label{sec:application_preliminaly}
The central tools of our argument are the centralizers $Z_W(W_I)$ and the normalizers $N_W(W_I)$ of standard parabolic subgroups $W_I$, which are described by the author \cite{Nui_centra} in a general setting (note that the normalizers had already been described by Brigitte Brink and Robert B.\ Howlett \cite{Bri-How}).
Here we summarize some of the author's results which we use.

Here we require the result only for the case that $|W_I|<\infty$ and $w_0(I) \in Z(W_I)$.
Now $Z_W(W_I)$ and $N_W(W_I)$ admit the following decompositions:
\begin{equation}
\label{eq:decompofZ}
Z_W(W_I)=(Z(W_I) \times W^{\perp I}) \rtimes Y_I \textrm{ and } N_W(W_I)=(W_I \times W^{\perp I}) \rtimes \widetilde{Y}_I.
\end{equation}
Here $W^{\perp I}$ denotes the subgroup of $W$ generated by the reflections in the set $Z_W(W_I) \smallsetminus W_I$, which is a Coxeter group by a theorem of Deodhar \cite{Deo_refsub} or of Dyer \cite{Dye}.
Since $Z(W_I)$ is an elementary abelian $2$-group, both $Z(W_I) \times W^{\perp I}$ and $W_I \times W^{\perp I}$ are also Coxeter groups.
The factor $\widetilde{Y}_I$ of $N_W(W_I)$ acts on $W_I \times W^{\perp I}$ as graph automorphisms, preserving the factor $W_I$.
The factor $Y_I$ of $Z_W(W_I)$ is torsion-free and is the kernel of the induced action of $\widetilde{Y}_I$ on $W_I$, so $Y_I$ is normal and has finite index in $\widetilde{Y}_I$ since $|W_I|<\infty$.
%%%%%
\subsubsection{The results}
\label{sec:application_statement}
Let $f:W \overset{\sim}{\to} W'$ be a group isomorphism between two Coxeter groups $W$ and $W'$ as above, and $I \subseteq S$ a subset with $|W_I|<\infty$ and $w_0(I) \in Z(W_I)$.
Our temporal subject is the element $f(w_0(I)) \in W'$.
Since $f(w_0(I))$ is an involution in $W'$ as well as $w_0(I)$, Theorem \ref{thm:Ric_involution} allows us to assume for simplicity that $f(w_0(I))=w_0(J)$ for some $J \subseteq S'$ with $|W'_J|<\infty$ and $w_0(J) \in Z(W'_J)$.
Let $Y'_J$ and $\widetilde{Y}'_J$ denote the last factors of $Z_{W'}(W'_J)$ and of $N_{W'}(W'_J)$, respectively (see (\ref{eq:decompofZ})).

We start with a very simple observation: since the isomorphism $f$ maps $w_0(I)$ to $w_0(J)$, it also maps $Z_W(w_0(I))$ onto $Z_{W'}(w_0(J))$, so the combination of Lemma \ref{lem:Zoflongestelement} and (\ref{eq:decompofZ}) yields the following isomorphism
\begin{equation}
\label{eq:relationbetweenZ}
f:(W_I \times W^{\perp I}) \rtimes \widetilde{Y}_I \overset{\sim}{\longrightarrow} (W'_J \times {W'}^{\perp J}) \rtimes \widetilde{Y}'_J.
\end{equation}
Let $\rho$ and $\rho'$ denote the maps representing the actions of $\widetilde{Y}_I$ and $\widetilde{Y}'_J$ in (\ref{eq:relationbetweenZ}), respectively.
Then by (\ref{eq:relationbetweenZ}) and the results in Section \ref{sec:application_preliminaly}, Theorem \ref{thm:cor_maintheorem} yields the following isomorphism
\begin{equation}
\label{eq:relationbetweenZ_2}
f:(W_I \times W^{\perp I})(\mathcal{O}_\rho) \rtimes (\widetilde{Y}_I)_\rho \overset{\sim}{\longrightarrow} (W'_J \times {W'}^{\perp J})(\mathcal{O}'_{\rho'}) \rtimes (\widetilde{Y}'_J)_{\rho'}.
\end{equation}
Now the left and the right sides of (\ref{eq:relationbetweenZ_2}) contain, as normal subgroups, $W_I$ and $W'_J$ which are $\rho(\widetilde{Y}_I)$-invariant and $\rho'(\widetilde{Y}'_J)$-invariant, respectively.
Thus if we know much enough of the structure of the left side of (\ref{eq:relationbetweenZ_2}), then we would be able to say something about the variation of the set $J$, so about the property of $f(w_0(I))$.
This is hopeful at least for individual cases, since \cite{Nui_centra} also gives a method for computing the explicit structure of the decompositions (\ref{eq:decompofZ}).

From now, we assume further that $\widetilde{Y}_I=Y_I$ (this is satisfied if $W_I$ admits no nontrivial graph automorphism).
For an arbitrary group $G$, let $G_{\mathrm{INV}}$ be the set of the involutions in $G$, so $\langle G_{\mathrm{INV}} \rangle \unlhd G$ and $\langle G_{\mathrm{INV}} \rangle$ is determined by the isomorphism type of $G$ only as well as $\langle G_{\mathrm{ACI}} \rangle$.
Then, since both $W_I \times W^{\perp I}$ and $W'_J \times {W'}^{\perp J}$ are generated by involutions and the torsion-free group $Y_I$ possesses no involution, we can derive from (\ref{eq:relationbetweenZ}) the following isomorphism
\begin{equation}
\label{eq:relationbetweenZ_3}
f:W_I \times W^{\perp I} \overset{\sim}{\longrightarrow} (W'_J \times {W'}^{\perp J}) \rtimes G,\quad \textrm{ where } G=\langle (\widetilde{Y}'_J)_{\mathrm{INV}} \rangle,
\end{equation}
by taking the $\langle (*)_{\mathrm{INV}} \rangle$ of both sides.
Now consider the centralizers of the normal subgroups $f^{-1}(W'_J)$ and $W'_J$ in the left and the right sides of (\ref{eq:relationbetweenZ_3}), respectively, which are also isomorphic via $f$.
Since $f^{-1}(W'_J)$ is generated by involutions, Proposition \ref{prop:Zofnormalinvolutions} implies that the centralizer in the left side is also generated by involutions, so is the centralizer in the right side.
The latter is the intersection of the right side of (\ref{eq:relationbetweenZ_3}) and $Z_{W'}(W'_J)=(Z(W'_J) \times {W'}^{\perp J}) \rtimes Y'_J$, that is
\[
(Z(W'_J) \times {W'}^{\perp J}) \rtimes (Y'_J \cap G),
\]
and all of its involutions are contained in the former factor since $Y'_J \cap G$ is torsion-free as well as $Y'_J$.
Thus it follows that $Y'_J \cap G=1$, so the $G$-action on the finite group $W'_J$ is faithful, therefore $G$ is also finite.
Hence, as mentioned in the first paragraph of Section \ref{sec:maintheorem_example}, (\ref{eq:relationbetweenZ_3}) and Theorem \ref{thm:cor_maintheorem} yield the following isomorphism
\begin{equation}
\label{eq:relationbetweenZ_4}
f:W_I \times W^{\perp I}{}_{\mathrm{fin}} \overset{\sim}{\longrightarrow} (W'_J \times {W'}^{\perp J}{}_{\mathrm{fin}}) \rtimes G_{\rho'}.
\end{equation}
This reduces our problem to the study of semidirect product decompositions of Coxeter groups whose irreducible components are finite.

Finally, specializing to the case $I=\{s\}$, we obtain the following result.
\begin{thm}
\label{thm:conditiontobereflection}
Let $(W,S)$ be an arbitrary Coxeter system.
\begin{enumerate}
\item Suppose that $s \in S$, and $W^{\perp s}{}_{\mathrm{fin}}$ is either trivial or generated by a single reflection conjugate to $s$.
Then $f(s) \in {S'}^{W'}$ for any Coxeter system $(W',S')$ and any group isomorphism $f:W \overset{\sim}{\to} W'$.
\item Suppose that every $s \in S$ satisfies the hypothesis of (1).
Then $f(S) \subseteq {S'}^{W'}$ for any Coxeter system $(W',S')$ and any group isomorphism $f:W \overset{\sim}{\to} W'$, so $f$ preserves the set of reflections.
Hence the set $S^W$ is determined by $W$ only and independent on the choice of $S \subseteq W$.
\end{enumerate}
\end{thm}
\begin{proof}
We only prove (1), since (2) follows immediately from (1) and the first remark of Section \ref{sec:maintheorem_application}.
Now the above argument works for $I=\{s\}$, so it suffices to deduce that $|J|=1$, implying that $f(s)=w_0(J) \in S'$ as desired.
This is immediately done if $W^{\perp s}{}_{\mathrm{fin}}=1$, since $J \neq \emptyset$ and now both sides of (\ref{eq:relationbetweenZ_4}) have cardinality $2$.

Suppose that $W^{\perp s}{}_{\mathrm{fin}}=\langle t \rangle$ with $t \in W$ conjugate to $s$.
Then both sides of (\ref{eq:relationbetweenZ_4}) have cardinality $4$.
Thus if $|J| \neq 1$, then it follows that $J=\{s',t'\}$ for two commuting generators $s',t' \in S'$ and the right side of (\ref{eq:relationbetweenZ_4}) is $W'_J$ itself, so we have an isomorphism $f:\langle s \rangle \times \langle t \rangle \overset{\sim}{\to} \langle s' \rangle \times \langle t' \rangle$.
Since we assumed that $f(s)=w_0(J)=s't'$, it follows that $f(t)$ is either $s'$ or $t'$, which cannot be conjugate to $f(s)=s't'$ in $W'$, contradicting the choice of $t$.
Hence $|J|=1$.
\end{proof}
Moreover, a forcecoming paper \cite{Nui_refindep} of the author will describe for which $s \in S$ the hypothesis is indeed satisfied, and show that this case occurs very frequently.
%%%%%
\section{Essential elements and Coxeter elements}
\label{sec:essentialelements}
Krammer introduced in his Ph.D.\ thesis \cite{Kra} the notion of essential elements of Coxeter groups.
An element $w$ of a Coxeter group $W$ is called \emph{essential} in $W$ if the parabolic closure $\mathrm{P}(w)$ of $w$ is $W$ itself (see Section \ref{sec:Coxetergroups_definition} for terminology).
Note that any $W$ of infinite rank cannot possess an essential element, while a Coxeter element $s_1s_2 \cdots s_n$ of an infinite irreducible $W$ of finite rank (where $S=\{s_1,s_2,\dots,s_n\}$) is always essential in $W$ (see Theorem \ref{thm:essentialelement}).
Here we summarize some properties of essential elements required in later sections, as follows:
\begin{thm}
\label{thm:essentialelement}
Let $W$ be an infinite irreducible Coxeter group of finite rank.
\begin{enumerate}
\item Any essential element of $W$ has infinite order.
\item Let $0 \neq k \in \mathbb{Z}$.
Then $w \in W$ is essential in $W$ if and only if $w^k$ is essential in $W$.
\item If $n=|S|$ and $\gamma_1,\dots,\gamma_n \in \Phi$ are linearly independent, then $s_{\gamma_1} \cdots s_{\gamma_n}$ is essential in $W$.
Hence any Coxeter element of $W$ is essential in $W$.
\end{enumerate}
\end{thm}
The claim (1) is an immediate consequence of a well-known theorem of Jacques Tits, which says that any finite subgroup of a Coxeter group is contained in a finite parabolic subgroup (see e.g.\ {\cite[Lemma 1.2]{Bah}} for a proof).
On the other hand, (2) and (3) are shown by Paris in his recent preprint \cite{Par}; however, he proved (3) only for Coxeter elements though his idea is adaptable applicable to the generalized version.
Here we include proofs of (2) and (3) along Paris' idea for the sake of completeness.

For (2), we fix $W$ and $w$ as in the statement.
For $\gamma \in \Phi$, let $\sigma^w_\gamma=((\sigma^w_\gamma)_n)_{n \in \mathbb{Z}}$ be the infinite sequence of $+$ and $-$ such that $(\sigma^w_\gamma)_n=\varepsilon$ if and only if $w^n \cdot \gamma \in \Phi^\varepsilon$.
We define $(\sigma^w_\gamma)_\infty$ (or $(\sigma^w_\gamma)_{-\infty}$, respectively) to be $\varepsilon \in \{+,-\}$ if $(\sigma^w_\gamma)_n=\varepsilon$ (or $(\sigma^w_\gamma)_{-n}=\varepsilon$, respectively) for all sufficiently large $n$.
Following \cite{Kra}, we say that $\gamma$ is \emph{$w$-periodic} if $w^n \cdot \gamma=\gamma$ for some $n \neq 0$.
Now we include the proofs of the following two lemmas for the sake of completeness.
\begin{lem}
[See {\cite[Proposition 5.2.2]{Kra}}]
\label{lem:charofoddroot}
If $\gamma \in \Phi$ is not $w$-periodic, then only finitely many sign-changes occur in the sequence $\sigma^w_\gamma$.
\end{lem}
\begin{proof}
By the hypothesis, all roots $w^n \cdot \gamma$ such that $(\sigma^w_\gamma)_n \neq (\sigma^w_\gamma)_{n+1}$ are distinct and contained in the finite set $\Phi\left[w\right] \cup -\Phi\left[w\right]$.
\end{proof}
A root $\gamma \in \Phi$ is called \emph{$w$-odd} (see \cite{Kra}) if it is not $w$-periodic (so both $(\sigma^w_\gamma)_{\pm \infty}$ are defined; see Lemma \ref{lem:charofoddroot}) and $(\sigma^w_\gamma)_\infty \neq (\sigma^w_\gamma)_{-\infty}$.
A reflection $s_\gamma$ is called \emph{$w$-odd} if $\gamma$ is $w$-odd.
\begin{lem}
[See {\cite[Lemma 5.2.7]{Kra}}]
\label{lem:oddrootforpower}
For $k \in \mathbb{Z} \smallsetminus \{0\}$, a root of $W$ is $w$-odd if and only if it is $w^k$-odd.
\end{lem}
\begin{proof}
Note that $\gamma \in \Phi$ is $w$-periodic if and only if it is $w^k$-periodic.
Thus for a non-$w$-periodic $\gamma$, all of $(\sigma^w_\gamma)_{\pm \infty}$ and $(\sigma^{w^k}_\gamma)_{\pm \infty}$ are defined (Lemma \ref{lem:charofoddroot}) and we have
\[
(\sigma^{w^k}_\gamma)_{\pm \infty}=
\begin{cases}
(\sigma^w_\gamma)_{\pm \infty} & \textrm{if } k>0;\\
(\sigma^w_\gamma)_{\mp \infty} & \textrm{if } k<0,
\end{cases}
\]
respectively.
Thus the claim follows.
\end{proof}
Let $\mathrm{P}^\infty(w)$ denote the subgroup of $W$ generated by the $w$-odd reflections.
The following result of \cite{Kra} is crucial in our argument.
\begin{prop}
[See {\cite[Corollary 5.8.7]{Kra}}]
\label{prop:parabolicclosureofw}
The parabolic closure $\mathrm{P}(w)$ is a direct product of $\mathrm{P}^\infty(w)$ and a finite number of finite groups.
\end{prop}
Moreover, the following result of the author \cite{Nui_indec} is also required.
See also {\cite[Theorem 4.1]{Par}} for the case of finite ranks.
\begin{prop}
[{\cite[Theorem 3.3]{Nui_indec}}]
\label{prop:indecomposability}
If $W$ is an infinite irreducible Coxeter group, then $W$ is directly indecomposable as an abstract group.
\end{prop}
\begin{cor}
\label{cor:charofessential}
Suppose that $W$ is infinite and irreducible.
Then $w \in W$ is essential in $W$ if and only if $\mathrm{P}^\infty(w)=W$.
\end{cor}
\begin{proof}
The `if' part is a consequence of Proposition \ref{prop:parabolicclosureofw}.
For the ``only if'' part, assume that $\mathrm{P}(w)=W$.
Then Proposition \ref{prop:parabolicclosureofw} implies that $W$ is the direct product of $\mathrm{P}^\infty(w)$ and certain finite groups, while $W$ is directly indecomposable (Proposition \ref{prop:indecomposability}).
Thus $W$ must coincide with one of the direct factors, which cannot be finite since $|W|=\infty$, so $W=\mathrm{P}^\infty(w)$ as desired.
\end{proof}
Now the claim (2) of Theorem \ref{thm:essentialelement} follows easily from Lemma \ref{lem:oddrootforpower} and Corollary \ref{cor:charofessential}, since the $w^k$-odd reflections are precisely the $w$-odd reflections.

For the proof of (3), we prepare two lemmas.
Here we say that \emph{$(W,S)$ is (non)degenerate} to signify the (non)degenerateness of the bilinear form $\langle \,,\, \rangle$, respectively.
\begin{lem}
[See {\cite[Lemma 3.2]{Par}}]
\label{lem:extendstonondegenerate}
Let $W$ be a Coxeter group of finite rank.
Then there is a nondegenerate Coxeter system $(\widetilde{W},\widetilde{S})$ of finite rank such that $S \subseteq \widetilde{S}$ and $\widetilde{W}_S=W_S$.
\end{lem}
\begin{proof}
We put $n=|S|$ and $S=\{s_1,s_3,\dots,s_{2n-1}\}$, and apply the following algorithm inductively for $1 \leq k \leq n$, beginning with $\widetilde{S}=S$:
\begin{quote}
if the Coxeter system $(\langle I_k \rangle,I_k)$ (where $I_k=\{s_i \in \widetilde{S} \mid i \leq 2k\}$) is degenerate, add a new generator $s_{2k}$ to $\widetilde{S}$ so that $s_{2k-1}s_{2k}$ has infinite order and $s_{2k}$ commutes with the other elements of $\widetilde{S}$.
\end{quote}
By computing the determinant of the matrix of the bilinear form with respect to the basis $\{\alpha_{s_i}\}_i$, it is checked inductively that the Coxeter system $(\langle I_k \rangle,I_k)$ will be nondegenerate when the $k$-th step is done.
Hence the Coxeter system $(\widetilde{W},\widetilde{S})=(\langle I_n \rangle,I_n)$ obtained finally is the desired one.
\end{proof}
\begin{lem}
\label{lem:1-eigenvector}
Any element of a proper standard parabolic subgroup $W_I$ of $W$ has a nonzero $1$-eigenvector in $V$.
\end{lem}
\begin{proof}
It suffices to consider the case that $S=\{s_1,s_2,\dots,s_n\}$ is finite and $I=S \smallsetminus \{s_n\}$.
Then, by definition of the $W$-action, the $n$-th row of the representation matrix $A_w$ of $w \in W_I$ relative to the basis $\Pi$ of $V$ is $(0\ 0\ \cdots\ 0\ 1)$.
Thus the matrix $I_n-A_w$ is singular as desired.
\end{proof}
The following property is the essence of the claim (3) of Theorem \ref{thm:essentialelement}.
\begin{prop}
\label{prop:independentrootforessential}
Let $W$ be a Coxeter group with $|S|=n<\infty$, and suppose that $\gamma_1,\dots,\gamma_n \in \Phi$ are linearly independent.
Then the standard parabolic closure of $s_{\gamma_1} \cdots s_{\gamma_n} \in W$ is $W$ itself.
\end{prop}
\begin{proof}
Assume contrary that $w=s_{\gamma_1} \cdots s_{\gamma_n} \in W_I$ for a proper $W_I \subset W$.
We may assume without loss of generality that $(W,S)$ is nondegenerate, since we can extend $S$ to $\widetilde{S}=S \sqcup \{t_1,\dots,t_m\}$ as in Lemma \ref{lem:extendstonondegenerate} and consider $t_1 \cdots t_mw \in \widetilde{W}_J$ instead of $w$, where $J=\widetilde{S} \smallsetminus (S \smallsetminus I)$.
Choose a nonzero $v \in V$ such that $w \cdot v=v$ (Lemma \ref{lem:1-eigenvector}).
Then, since $(W,S)$ is nondegenerate, there is an index $i$ such that $\langle v,\gamma_i \rangle \neq 0$ and $\langle v,\gamma_j \rangle=0$ for all $j>i$.
This implies that $w \cdot v=s_{\gamma_1} \cdots s_{\gamma_i} \cdot v$, which is the sum of $s_{\gamma_i} \cdot v=v-2\langle v,\gamma_i \rangle \gamma_i$ and a linear combination of $\gamma_1,\dots,\gamma_{i-1}$.
Now the property $w \cdot v=v$ yields an expression of $2\langle v,\gamma_i \rangle \gamma_i$ as a linear combination of the other $\gamma_j$, contradicting the linear independence of $\gamma_1,\dots,\gamma_n$.
Hence the claim follows.
\end{proof}
Now the claim (3) of Theorem \ref{thm:essentialelement} is easily proved, since the hypothesis of Proposition \ref{prop:independentrootforessential} is invariant under the action of $W$.
Hence the proof of Theorem \ref{thm:essentialelement} is concluded.
%%%%%
\section{On the fixed-point subgroups by Coxeter graph automorphisms}
\label{sec:fixedpoint}
The subject of this section is the fixed-point subgroup
\[
W^\tau=\{w \in W \mid \tau(w)=w\}
\]
of a Coxeter group $W$ by a graph automorphism $\tau \in \mathrm{Aut}\,\Gamma$ (as mentioned in Section \ref{sec:Coxetergroups_definition}, the automorphism of $W$ induced by $\tau$ is also denoted by $\tau$).
Let $\tau \backslash S$ denote the set of the $\langle \tau \rangle$-orbits in $S$.
Then it was shown by Steinberg {\cite[Theorem 1.32]{Ste}} that $W^\tau$ is a Coxeter group with respect to the following generating set
\[
S(W^\tau)=\{w_0(I) \in W \mid I \in \tau \backslash S \textrm{ and } |W_I|<\infty\}
\]
(see also \cite{Muh} and \cite{Nan}).
Here we show the following properties of the subgroup $W^\tau$, which will be used in the proof of the main theorem.
\begin{thm}
\label{thm:indexoffixedpointsubgroup}
Let $W$ be an arbitrary Coxeter group and $\tau \in \mathrm{Aut}\,\Gamma$.
Then $W^\tau$ has finite index in $W$ if and only if $\tau$ is identity on all irreducible components of $W$ except a finite number of finite irreducible components.
\end{thm}
\begin{thm}
\label{thm:Wtauisinfinite}
Let $W$ be an infinite irreducible Coxeter group and $\tau \in \mathrm{Aut}\,\Gamma$.
\begin{enumerate}
\item If $|W_I|<\infty$ for all $I \in \tau \backslash S$, then the Coxeter group $W^\tau$ is also infinite and irreducible with respect to the generating set $S(W^\tau)$.
\item Suppose that the hypothesis of (1) fails and every orbit $I \in \tau \backslash S$ is finite.
Then for any $1 \neq w \in W^\tau$, there is an element $u \in W$ of infinite order such that $u^kw\tau(u)^{-k} \neq w$ for all $0 \neq k \in \mathbb{Z}$.
\end{enumerate}
\end{thm}
Note that the result on infiniteness of $W^\tau$ in Theorem \ref{thm:Wtauisinfinite} (1) is mentioned in {\cite[Section 5]{Muh}} without proof in a generalized setting.
%%%%%
\subsection{Proof of Theorem \ref{thm:indexoffixedpointsubgroup}}
\label{sec:proof_indexoffixedpointsubgroup}
Our first step is to prove the following lemma:
\begin{lem}
\label{lem:specialcase_twoinfiniteorder}
Let $W$ be an (irreducible) affine or compact hyperbolic Coxeter group with $\mathrm{type}\,W \neq \widetilde{A_1}$ (see Section \ref{sec:hyperbolicCoxetergroups} for terminology).
Suppose further that $\mathrm{Aut}\,\Gamma \neq \{\mathrm{id}_S\}$.
Then for any $\mathrm{id}_S \neq \tau \in \mathrm{Aut}\,\Gamma$, there is an element $w \in W$ of infinite order such that $\langle w \rangle \cap \langle \tau(w) \rangle=1$.
\end{lem}
From now until the end of the proof of Lemma \ref{lem:specialcase_twoinfiniteorder}, we assume that $S$ is finite and the base field of the (finite-dimensional) geometric representation space $V$ is extended from $\mathbb{R}$ to $\mathbb{C}$.
Then the bilinear form $\langle \,,\, \rangle$ and the faithful $W$-action also extend naturally so that $W$ is embedded injectively in the group of orthogonal linear transformations of $V$ relative to $\langle \,,\, \rangle$.
For $\lambda \in \mathbb{C}$, let $V_\lambda(w)$ denote the $\lambda$-eigenspace of $w \in W$, and let $V_{\sqrt{1}}(w)$, $V_{\neq \sqrt{1}}(w)$ be the sum of $V_\lambda(w)$ where $\lambda$ runs over the roots of unity, over $\mathbb{C} \smallsetminus \{0\}$ except the roots of unity, respectively.
Then some elementary linear algebra shows that, if $w \in W$, $0 \neq \lambda \in \mathbb{C}$ and $0 \neq k \in \mathbb{Z}$, then $V_\lambda(w^k)$ is the sum of $V_\mu(w)$ where $\mu \in \mathbb{C}$ varies subject to $\mu^k=\lambda$.
Hence we have $V_{\sqrt{1}}(w^k)=V_{\sqrt{1}}(w)$ and $V_{\neq \sqrt{1}}(w^k)=V_{\neq \sqrt{1}}(w)$ whenever $k \neq 0$.

Now we have the following:
\begin{lem}
\label{lem:eigenspaceandintersection}
Let $w_1,w_2 \in W$ and suppose that either $V_{\sqrt{1}}(w_1) \neq V_{\sqrt{1}}(w_2)$ or $V_{\neq \sqrt{1}}(w_1) \neq V_{\neq \sqrt{1}}(w_2)$.
Then $\langle w_1 \rangle \cap \langle w_2 \rangle=1$.
\end{lem}
\begin{proof}
Assume contrary that $k,\ell \in \mathbb{Z} \smallsetminus \{0\}$ and $w_1{}^k=w_2{}^\ell$.
Then, in the first case $V_{\sqrt{1}}(w_1) \neq V_{\sqrt{1}}(w_2)$, the above observation implies that
\[
V_{\sqrt{1}}(w_1{}^k)=V_{\sqrt{1}}(w_1) \neq V_{\sqrt{1}}(w_2)=V_{\sqrt{1}}(w_2{}^\ell),
\]
contradicting the assumption $w_1{}^k=w_2{}^\ell$.
The other case is similar.
\end{proof}
Define actions of $\tau \in \mathrm{Aut}\,\Gamma$ on $V$ and the dual space $V^*$ with dual basis $\{\alpha_s^* \mid s \in S\}$ (as linear transformations) by
\[
\tau(\alpha_s)=\alpha_{\tau(s)} \textrm{ and } \tau(\alpha_s^*)=\alpha_{\tau(s)}^* \textrm{ for } s \in S.
\]
Then $\tau$ preserves the bilinear form $\langle \,,\, \rangle$, and we have $\tau(w) \cdot \tau(v)=\tau(w \cdot v)$ for $w \in W$ and $v \in V$.
Thus for $0 \neq \lambda \in \mathbb{C}$ and $w \in W$, it follows that $V_\lambda(\tau(w))=\tau(V_\lambda(w))$, $V_{\sqrt{1}}(\tau(w))=\tau(V_{\sqrt{1}}(w))$ and $V_{\neq \sqrt{1}}(\tau(w))=\tau(V_{\neq \sqrt{1}}(w))$.
Moreover, we have $\tau(\eta)(\tau(v))=\eta(v)$ for $\eta \in V^*$ and $v \in V$.
Note also that $\mathrm{Ann}(\tau(V'))=\tau(\mathrm{Ann}(V'))$ for any subspace $V' \subseteq V$, where $\mathrm{Ann}(V')=\{\eta \in V^* \mid \eta(V')=0\}$ denotes the annihilator of $V'$.\\
\indent
By these observations, we have the following lemmas.
In these lemmas, write $v^\perp=\{v' \in V \mid \langle v,v' \rangle=0\}$ for $v \in V$.
\begin{lem}
\label{lem:graphautoandinfiniteorder}
Let $\mathrm{id}_S \neq \tau \in \mathrm{Aut}\,\Gamma$, $\beta,\gamma \in \Phi^+$ and $V' \subset V$ a subspace of codimension $1$.
Suppose that $\langle \beta,\gamma \rangle=-1$, $V' \subseteq \beta^\perp \cap \gamma^\perp$ and $\mathrm{Ann}(V')$ is not $\tau$-invariant.
Then $w=s_\beta s_\gamma \in W$ has infinite order and $\langle w \rangle \cap \langle \tau(w) \rangle=1$.
\end{lem}
\begin{proof}
Since $\langle \beta,\gamma \rangle=-1$, we have $w^k \cdot \beta=(2k+1)\beta+2k\gamma \neq \beta$ for all $k \geq 1$, showing that $w$ has infinite order.
Thus $V_{\sqrt{1}}(w) \neq V$, since otherwise we have $V_1(w^k)=V$ and $w^k=1$ for a sufficiently large $k$, a contradiction.
Now we have
\[
V' \subseteq \beta^\perp \cap \gamma^\perp \subseteq V_1(w) \subseteq V_{\sqrt{1}}(w) \subset V \textrm{ and } \mathrm{dim}\,V-\mathrm{dim}\,V'=1,
\]
implying that $V'=V_{\sqrt{1}}(w)$.
Since $\mathrm{Ann}(V')$ is not $\tau$-invariant, we have
\[
\mathrm{Ann}(V_{\sqrt{1}}(w)) \neq \tau(\mathrm{Ann}(V_{\sqrt{1}}(w)))=\mathrm{Ann}(V_{\sqrt{1}}(\tau(w))),
\]
so $V_{\sqrt{1}}(w) \neq V_{\sqrt{1}}(\tau(w))$.
Hence Lemma \ref{lem:eigenspaceandintersection} completes the proof.
\end{proof}
\begin{lem}
\label{lem:lemmaforinfiniteorder}
For $i=1,2$, let $\beta_i,\gamma_i \in \Phi^+$ and $V^{(i)} \subset V$ a subspace of codimension $3$, and suppose that $\langle \beta_i,\gamma_i \rangle<-1$, $V^{(i)} \subseteq {\beta_i}^\perp \cap {\gamma_i}^\perp$ and $\mathbb{C}\beta_1+\mathbb{C}\gamma_1 \neq \mathbb{C}\beta_2+\mathbb{C}\gamma_2$.
Then each $w_i=s_{\beta_i}s_{\gamma_i} \in W$ has infinite order and $\langle w_1 \rangle \cap \langle w_2 \rangle=1$.
\end{lem}
\begin{proof}
Put ${v_i}^\pm=(-c_i \pm \sqrt{c_i^2-1})\beta_i+\gamma_i$ and ${\lambda_i}^\pm=2c_i^2-1 \mp 2c_i\sqrt{c_i^2-1}$, respectively, where $c_i=\langle \beta_i,\gamma_i \rangle$.
Then a direct computation shows that $w_i \cdot {v_i}^\pm={\lambda_i}^\pm {v_i}^\pm$ and $|{\lambda_i}^\pm| \neq 1$, respectively, and ${\lambda_i}^+{\lambda_i}^-=1$, so $w_i$ has infinite order.
Moreover, since ${\beta_i}^\perp \cap {\gamma_i}^\perp \subseteq V_1(w_i)$, the hypothesis implies that $\mathrm{dim}\,V-\mathrm{dim}\,V_1(w_i) \leq 3$, so the characteristic polynomial $\chi_{w_i}(x)=\mathrm{det}(x \cdot \mathrm{id}_V-w_i)$ of $w_i$ decomposes as
\[
\chi_{w_i}(x)=(x-1)^{|S|-3}(x-{\lambda_i}^+)(x-{\lambda_i}^-)(x-\mu_i) \textrm{ where } \mu_i \in \mathbb{C}.
\]
Now we have $\pm 1=\mathrm{det}\,w_i=\pm \chi_{w_i}(0)=\pm {\lambda_i}^+{\lambda_i}^-\mu_i$ since $w_i$ is a product of involutions, so $\mu_i=\pm 1$.
Thus $V_{\neq \sqrt{1}}(w_i)=\mathbb{C}{v_i}^++\mathbb{C}{v_i}^-=\mathbb{C}\beta_i+\mathbb{C}\gamma_i$, so $V_{\neq \sqrt{1}}(w_1) \neq V_{\neq \sqrt{1}}(w_2)$ by the hypothesis.
Hence Lemma \ref{lem:eigenspaceandintersection} completes the proof.
\end{proof}
\begin{cor}
\label{cor:graphautoandinfiniteorder_2}
Let $\mathrm{id}_S \neq \tau \in \mathrm{Aut}\,\Gamma$, $\beta,\gamma \in \Phi^+$ and $V' \subset V$ a subspace of codimension $3$.
Suppose that $\langle \beta,\gamma \rangle<-1$, $V' \subseteq \beta^\perp \cap \gamma^\perp$ and $\mathbb{C}\beta+\mathbb{C}\gamma$ is not $\tau$-invariant.
Then $w=s_\beta s_\gamma \in W$ has infinite order and $\langle w \rangle \cap \langle \tau(w) \rangle=1$.
\end{cor}
\begin{proof}
Note that $\tau(\mathbb{C}\beta+\mathbb{C}\gamma)=\mathbb{C}\tau(\beta)+\mathbb{C}\tau(\gamma)$ and $\tau(s_\beta s_\gamma)=s_{\tau(\beta)}s_{\tau(\gamma)}$.
Then the claim follows from Lemma \ref{lem:lemmaforinfiniteorder}, where $\beta_1=\beta$, $\gamma_1=\gamma$, $\beta_2=\tau(\beta)$, $\gamma_2=\tau(\gamma)$, $V^{(1)}=V'$ and $V^{(2)}=\tau(V')$.
\end{proof}

\noindent
\textbf{Proof of Lemma \ref{lem:specialcase_twoinfiniteorder}.}
This lemma is deduced from Lemma \ref{lem:graphautoandinfiniteorder} for affine case and Corollary \ref{cor:graphautoandinfiniteorder_2} for compact hyperbolic case, by constructing the $\beta,\gamma$ and $V'$ as in Tables \ref{tab:list_twoinfiniteorder_affine} and \ref{tab:list_twoinfiniteorder_hyperbolic} (see also Figures \ref{fig:affine} and \ref{fig:compacthyperbolic}).
Note that $\beta+\gamma$ is the null root of $W$ in an affine case.
If $|\mathrm{Aut}\,\Gamma| \geq 3$, we assume by symmetry that $\tau$ satisfies the condition in the second column of the lists, where we abbreviate $s_i$ to $i$.
In the next two columns, a word $c_1c_2 \cdots c_r$ (where $r=|S|$) signifies $\sum_{i=1}^{r}c_i\alpha_i \in V$ and $\widetilde{\alpha}_i$ denotes the unique highest root of the finite Coxeter group $W_{S \smallsetminus \{s_i\}}$.
Finally, the last column gives a basis of $V'$ or of $\mathrm{Ann}(V')$.\qed
%%
%%%%%
\begin{table}
\centering
\caption{List for the proof of Lemma \ref{lem:specialcase_twoinfiniteorder}, affine case}
\label{tab:list_twoinfiniteorder_affine}\medskip
\begin{tabular}{c|c||c|c|c}
$W$ & $\tau$ & $\beta$ & $\gamma$ & $\mathrm{Ann}(V')$ \\ \hline
\vbox to1.2em{}$\widetilde{A_n}$ ($n \geq 2$) & $\tau(1) \neq 1$ & $\alpha_1$ & $011 \cdots 11$ & $2\alpha_1^*-\alpha_2^*-\alpha_{n+1}^*$ \\
$\widetilde{B_n}$ ($n \geq 3$) & & $\alpha_1$ & $0122 \cdots 22\sqrt{2}$ & $2\alpha_1^*-\alpha_3^*$ \\
$\widetilde{C_n}$ ($n \geq 2$) & & $\alpha_1$ & $0\sqrt{2}\sqrt{2} \cdots \sqrt{2}\sqrt{2}1$ & $\sqrt{2}\alpha_1^*-\alpha_2^*$ \\
$\widetilde{D_n}$ ($n \geq 4$) & $\tau(1) \neq 1$ & $\alpha_1$ & $0122 \cdots 2211$ & $2\alpha_1^*-\alpha_3^*$ \\
$\widetilde{E_6}$ & $\tau(1) \neq 1$ & $\alpha_1$ & $0232121$ & $2\alpha_1^*-\alpha_2^*$ \\
$\widetilde{E_7}$ & & $\alpha_1$ & $02343212$ & $2\alpha_1^*-\alpha_2^*$
\end{tabular}
\end{table}
%%%%%
\begin{table}
\centering
\caption{List for the proof of Lemma \ref{lem:specialcase_twoinfiniteorder}, compact hyperbolic case}
\label{tab:list_twoinfiniteorder_hyperbolic}\medskip
\begin{tabular}{c|c||c|c|c}
$W$ & $\tau$ & $\beta$ & $\gamma$ & $V'$ \\ \hline
\vbox to1.2em{}$X_1$ & & $\alpha_1$ & $01111$ & $\alpha_3,\alpha_4$ \\
$X_2(m_1,m_2)$ & & $\alpha_4$ & $s_1s_2 \cdot \alpha_3$ & $\alpha_2$ \\
$X_3(m_1,m_2,m_3)$ & $\tau(3) \neq 3$ & $\alpha_3$ & $s_2 \cdot \alpha_1$ & $\emptyset$ \\
$Y_1$ & & $\alpha_4$ & $\widetilde{\alpha}_4$ & $\alpha_3$ \\
$Y_2$ & & $\alpha_5$ & $\widetilde{\alpha}_5$ & $\alpha_1,\alpha_2$ \\
$Y_3(5)$ & & $\alpha_5$ & $\widetilde{\alpha}_5$ & $\alpha_1,\alpha_2$ \\
$Y_4(5)$ & & $\alpha_4$ & $\widetilde{\alpha}_4$ & $\alpha_1$ \\
$Y_5$ & & $\alpha_4$ & $\widetilde{\alpha}_4$ & $\alpha_1$ \\
$Y_6(m,m)$ ($m \geq 5$) & & $\alpha_3$ & $s_2 \cdot \alpha_1$ & $\emptyset$
\end{tabular}
\end{table}
%%%%%

Now we cancel the assumption $|S|<\infty$ placed above.
To prove Theorem \ref{thm:indexoffixedpointsubgroup}, note that if $\tau \in \mathrm{Aut}\,\Gamma$ leaves $W_I \subseteq W$ invariant, then $W_I$ possesses its own fixed-point subgroup ${W_I}^\tau$ which coincides with $W^\tau \cap W_I$.\medskip

\noindent
\textbf{Proof of Theorem \ref{thm:indexoffixedpointsubgroup}.}
The only nontrivial part is the ``only if'' part, so we prove it.
Note that, by (\ref{eq:index_intersection}), the hypothesis implies that
\begin{equation}
\label{eq:proof_indexoffixedpointsubgroup}
\left[G:W^\tau \cap G\right]<\infty \textrm{ for any subgroup } G \leq W,
\end{equation}
so $W^\tau \cap G \neq 1$ for every infinite subgroup $G$ of $W$.\medskip\\
\textbf{Step 1: if $I \subseteq S$ is finite, and $W_I$ is infinite and irreducible, then $\tau(W_I)=W_I$.}\medskip

Assume contrary that $\tau(I) \neq I$, or equivalently $I \not\subseteq \tau(I)$.
Then we have $I \cap \tau(I) \neq I$, while $W^\tau \cap W_I \subseteq W_I \cap W_{\tau(I)}=W_{I \cap \tau(I)}$ (see (\ref{eq:intersectionofW_I})), therefore no essential element in $W_I$ lies in $W^\tau$.
Hence by Theorem \ref{thm:essentialelement}, any power $w^k$ (with $k \neq 0$) of a Coxeter element $w$ of $W_I$ has infinite order and is not in $W^\tau$, so we have $W^\tau \cap \langle w \rangle=1$, contradicting (\ref{eq:proof_indexoffixedpointsubgroup}).\medskip\\
\textbf{Step 2: the claim holds if $W$ has finite rank.}\medskip

Now it suffices to show that $\tau$ is identity on every infinite irreducible component $W_I$.
Moreover, since (by Step 1) $\tau(W_I)=W_I$ and (by (\ref{eq:proof_indexoffixedpointsubgroup})) $\left[W_I:{W_I}^\tau\right]<\infty$, it actually suffices to consider the case $W_I=W$, namely $W$ itself is infinite and irreducible.
In this case, our aim is to show that $\tau$ is identity.

First, we consider the case that $W$ is not of type $\widetilde{A_1}$ and every proper $W_J \subset W$ is infinite.
Then by combining Proposition \ref{prop:charofminimalnonfinite} (2) and Lemma \ref{lem:specialcase_twoinfiniteorder}, we have $\langle w \rangle \cap \langle \tau(w) \rangle=1$ for some $w \in W$ of infinite order whenever $\tau \neq \mathrm{id}_S$.
This implies that $W^\tau \cap \langle w \rangle=1$, contradicting (\ref{eq:proof_indexoffixedpointsubgroup}).
Thus $\tau$ must be identity now, as desired.
On the other hand, the claim also holds if $\mathrm{type}\,W=\widetilde{A_1}$, since now we have $W^\tau=1$ whenever $\tau \neq \mathrm{id}_S$.

Finally, we consider the remaining case that a proper $W_J \subset W$ is infinite.
We may assume that $J=S \smallsetminus \{s\}$ for some $s \in S$, so it suffices to show that $\tau|_J=\mathrm{id}_J$.
Since $|S|<\infty$, we may assume further that $W_J$ is irreducible: indeed, if $W_J$ is not irreducible and $W_K$ is an infinite irreducible component of $W_J$ (which exists since $|J|<\infty$), then the set $S \smallsetminus \{s'\}$, where $s'$ is an element of $J \smallsetminus K$ farthest from $s$ in $\Gamma$, possesses the desired properties.
Now Step 1 implies that $\tau(J)=J$, so $\left[W_J:{W_J}^\tau\right]<\infty$ (by (\ref{eq:proof_indexoffixedpointsubgroup})), therefore the induction on $|S|$ shows that $\tau|_J=\mathrm{id}_J$, as desired.\medskip\\
\textbf{Step 3: if $I \in \tau \backslash S$, then $|I|<\infty$.}\medskip

Assume contrary that $|I|=\infty$.
Then for any $w \in W_I$ with $J=\mathrm{supp}(w)$ (finite and) nonempty, we have $J \neq I$ and so $J \neq \tau(J)$ (since $I$ is a $\langle \tau \rangle$-orbit), therefore $J \not\subseteq \tau(J)$ and $w \not\in W_{\tau(J)}$.
This means that $\tau(w) \neq w$.
Thus we have $W^\tau \cap W_I=1$, contradicting (\ref{eq:proof_indexoffixedpointsubgroup}).\medskip\\
\textbf{Step 4: $\tau$ is identity on every infinite irreducible component $W_I$.}\medskip

First, we consider the case that a (not necessarily proper) $W_J \subseteq W_I$ of finite rank is infinite.
We can take an irreducible $W_J$.
Now assume contrary that $\tau$ is not identity on $W_I$, so $\tau(s) \neq s$ for some $s \in I$.
Then, since $W_I$ is irreducible and $|J|<\infty$, an irreducible $W_K \subseteq W_I$ of finite rank contains both $W_J$ and $s$.
This $W_K$ is also infinite, so $\tau(K)=K$ (Step 1), therefore $\left[W_K:{W_K}^\tau\right]<\infty$ by (\ref{eq:proof_indexoffixedpointsubgroup}).
Now Step 2 implies that $\tau$ is identity on $W_K$, contradicting the choice of $s$.
Hence the claim holds in this case.

In the remaining case, $W_I$ is of type $A_\infty$, $A_{\pm \infty}$, $B_\infty$ or $D_\infty$ (see Figure \ref{fig:nonf.g.Coxetergroups}) as mentioned in Section \ref{sec:hyperbolicCoxetergroups}.
Note that $\tau(W_I)=W_I$, since otherwise we have $W^\tau \cap W_I=1$, contradicting (\ref{eq:proof_indexoffixedpointsubgroup}).
Now the claim is trivial in the first and the third cases where $\mathrm{Aut}\,\Gamma=\{\mathrm{id}_S\}$.

In the case $\mathrm{type}\,W_I=A_{\pm \infty}$, if $\tau$ is not identity on $W_I$, then Step 3 implies that $\tau$ is a turning of the infinite path $\Gamma_I$, so there is an infinite $J \subset I$ with $J \cap \tau(J)=\emptyset$.
Now we have $W^\tau \cap W_J=1$, contradicting (\ref{eq:proof_indexoffixedpointsubgroup}).
This verifies the claim.

Finally, in the case $\mathrm{type}\,W_I=D_\infty$, if $\tau$ is not identity on $W_I$, then $\tau(s_1)=s_2$, $\tau(s_2)=s_1$ and $\tau$ fixes $J=I \smallsetminus \{s_1,s_2\}$ pointwise.
Put $K=J \cup \{s_2\}$.
Since any $w \in W_K$ satisfies that $\tau(w) \in W_{J \cup \{s_1\}}$, we have $W^\tau \cap W_K=W_J$ (see (\ref{eq:intersectionofW_I})), so $\left[W_K:W_J\right]<\infty$ by (\ref{eq:proof_indexoffixedpointsubgroup}).
However, putting $\gamma_k=\sum_{i=2}^{k}\alpha_{s_i} \in \Phi_K^+$ for $k \geq 3$, Lemma \ref{lem:reflectionincoset} implies that all of the infinitely many reflections $s_{\gamma_k}$ belong to distinct cosets in $W_K/W_J$.
This contradiction yields the claim.\medskip\\
\textbf{Step 5: conclusion.}\medskip

Assume that the ``only if'' part fails.
Then by Step 4, $W$ possesses infinitely many finite irreducible components $W_{I_1},W_{I_2},\dots$ on which $\tau$ is not identity.
Since every $\langle \tau \rangle$-orbit is finite (Step 3), there is an infinite sequence $s_1,s_2,\dots$ of distinct elements of $S$ such that $J=\{s_i \mid i \geq 1\}$ satisfies that $\tau(J) \cap J=\emptyset$; take $s_1$ as any element of $I_1$ with $\tau(s_1) \neq s_1$, and if $s_1,\dots,s_k$ are already chosen, then take $s_{k+1} \in I_i$ where $I_i$ does not intersect with the $\langle \tau \rangle$-orbits of the preceding $s_j$ and $\tau(s_{k+1}) \neq s_{k+1}$.
Now we have $W^\tau \cap W_J=1$ and $|W_J|=\infty$, contradicting (\ref{eq:proof_indexoffixedpointsubgroup}).
Hence the proof is concluded.\qed
%%
%%%%%
\subsection{Proof of Theorem \ref{thm:Wtauisinfinite}}
\label{sec:proof_Wtauisinfinite}
We start with some preliminaries.
Let $\tau \in \mathrm{Aut}\,\Gamma$ and $w \in W^\tau$, and denote the support of $w$ as an element of $(W^\tau,S(W^\tau))$ by $\mathrm{supp}^\tau(w)$.
The following (part of a) result of \cite{Nan} shows a relation between $\mathrm{supp}(w)$ and $\mathrm{supp}^\tau(w)$.
\begin{prop}
[See {\cite[Proposition 3.3]{Nan}}]
\label{prop:relationofsupport}
Let $w=w_0(I_1) \cdots w_0(I_r)$ (where $w_0(I_i) \in S(W^\tau)$) be a reduced expression of $w \in W^\tau$ with respect to $S(W^\tau)$.
Then any expression of $w$ obtained by replacing each $w_0(I_i)$ with its reduced expression, with respect to $S$, is also reduced with respect to $S$.
Hence $\mathrm{supp}(w)=\bigcup_{i=1}^rI_i$.
\end{prop}
Secondly, we give a remark on the Coxeter graph of the Coxeter system $(W^\tau,S(W^\tau))$, denoted here by $\Gamma^\tau$.
Let $\tau \backslash \Gamma$ be the graph with vertex set $\tau \backslash S$, in which two orbits $I,J \in \tau \backslash S$ are joined if and only if these sets are adjacent in $\Gamma$.
Then the vertex set $S(W^\tau)$ of $\Gamma^\tau$ is regarded as a subset of the vertex set $\tau \backslash S$ of $\tau \backslash \Gamma$ via an embedding $w_0(I) \mapsto I$.
Now we have the following result on a relation between $\Gamma^\tau$ and $\tau \backslash \Gamma$.
\begin{lem}
\label{lem:Gammatauisfullsubgraph}
Under the embedding $S(W^\tau) \hookrightarrow \tau \backslash S$ of the vertex set, the underlying graph of $\Gamma^\tau$ is a full subgraph of $\tau \backslash \Gamma$.
\end{lem}
\begin{proof}
Let $I,J \in \tau \backslash S$ be two distinct orbits with both $W_I$ and $W_J$ finite.
It is obvious that $I$ and $J$ are not adjacent in $\Gamma^\tau$ (i.e.\ $w_0(I)$ and $w_0(J)$ commute) if these are not adjacent in $\tau \backslash S$.
Thus our remaining task is to show that $w_0(I)$ and $w_0(J)$ do not commute if $I$ and $J$ are adjacent in $\tau \backslash S$, namely some $s \in I$ is adjacent to $J$ in $\Gamma$.
Now Lemma \ref{lem:lemmaforsupport} (1) implies that $w_0(J) \cdot \alpha_s \in \Phi_{J \cup \{s\}}^+ \smallsetminus \Phi_{\{s\}}$, so $w_0(I)w_0(J) \cdot \alpha_s \in \Phi^+$.
On the other hand, we have $w_0(I) \cdot \alpha_s \in \Phi_I^-$, so $w_0(J)w_0(I) \cdot \alpha_s \in \Phi^-$.
Thus we have $w_0(I)w_0(J) \neq w_0(J)w_0(I)$ as desired.
\end{proof}
Moreover, note that $\tau \backslash \Gamma$ is connected whenever $\Gamma$ is.
Indeed, for any $I,J \in \tau \backslash S$, a path in the connected graph $\Gamma$ between any $s \in I$ and any $t \in J$ gives rise to a path in $\tau \backslash \Gamma$ between $I$ and $J$.\medskip

\noindent
\textbf{Proof of Theorem \ref{thm:Wtauisinfinite} (1).}
As is remarked above, the irreducibility of $W$ yields the connectedness of $\tau \backslash \Gamma$, while the hypothesis implies that the embedding $\Gamma^\tau \hookrightarrow \tau \backslash \Gamma$ in Lemma \ref{lem:Gammatauisfullsubgraph} is now an isomorphism.
Thus $\Gamma^\tau$ is connected as desired.

For the infiniteness of $W^\tau$, assume the contrary.
Then $W^\tau$ possesses the longest element $w_0^\tau$ with respect to $S(W^\tau)$.
Now for any $s \in S$, belonging by the hypothesis to an $I \in \tau \backslash S$ with $|W_I|<\infty$, the $w_0^\tau$ and $w_0(I)$ admit a reduced expression with respect to $S(W^\tau)$ and $S$ ending with $w_0(I)$ and $s$, respectively, by Exchange Condition.
Thus Proposition \ref{prop:relationofsupport} implies that $w_0^\tau$ admits a reduced expression with respect to $S$ ending with $s$.
Since $s \in S$ is arbitrary, this means that $W$ is finite and $w_0^\tau$ is the longest element of $W$ (see Section \ref{sec:hyperbolicCoxetergroups}), contradicting the hypothesis that $W$ is infinite.
Hence the claim follows.\qed\medskip

\noindent
\textbf{Proof of Theorem \ref{thm:Wtauisinfinite} (2).}
Note that the graph $\tau \backslash \Gamma$ is connected.
Since the hypothesis of (1) now fails, there is a path $I_0I_1 \cdots I_r$ in $\tau \backslash \Gamma$, where $I_i \in \tau \backslash S$, such that $w_0(I_0) \in \mathrm{supp}^\tau(w)$ and $|W_{I_r}|=\infty$.
By choosing the shortest possible path, we may assume that $|W_{I_i}|<\infty$ and $w_0(I_i) \not\in \mathrm{supp}^\tau(w)$ for $1 \leq i \leq r-1$.
Now Lemma \ref{lem:Gammatauisfullsubgraph} says that $I_0I_1 \cdots I_{r-1}$ is also a path in $\Gamma^\tau$, so by applying Lemma \ref{lem:lemmaforsupport} (2) to the Coxeter system $(W^\tau,S(W^\tau))$, it is deduced that $w_0(I_{r-1}) \in \mathrm{supp}^\tau(w'w{w'}^{-1})$ where $w'=w_0(I_{r-1}) \cdots w_0(I_2)w_0(I_1) \in W^\tau$.
Thus Proposition \ref{prop:relationofsupport} implies that $\mathrm{supp}(w'w{w'}^{-1}) \subseteq S$ contains $I_{r-1}$, does not intersect $I_r$ and is adjacent to $I_r$ in $\Gamma$.

Take $s \in I_r$ adjacent to $\mathrm{supp}(w'w{w'}^{-1})$ in $\Gamma$.
Now we show that, if $W_{I_r}$ possesses an element $u'$ of infinite order such that $s \in \mathrm{supp}({u'}^k)$ for all $0 \neq k \in \mathbb{Z}$, then $u={w'}^{-1}\tau^{-1}(u')w'$ is the desired element.
Indeed, for $k \neq 0$, we have $\tau^{-1}(u')^kw'w{w'}^{-1}{u'}^{-k} \neq w'w{w'}^{-1}$ by the choice of $u'$ and Lemma \ref{lem:lemmaforsupport} (3) (note that $\tau^{-1}(u') \in W_{I_r}$), so, since $\tau(w')=w'$, we have
\[
u^kw\tau(u)^{-k}={w'}^{-1}\bigl(\tau^{-1}(u')^kw'w{w'}^{-1}{u'}^{-k}\bigr)w' \neq w.
\]

Finally, we show the existence of such an element $u'$, concluding the proof.
Since $|W_{I_r}|=\infty$ and $I_r \in \tau \backslash S$ is a finite orbit, an irreducible component of $W_{I_r}$, therefore that containing $s$, is infinite.
Now Theorem \ref{thm:essentialelement} implies that a Coxeter element of this component possesses the desired property.\qed
%%
%%%%%
\section{Proof of the main theorem}
\label{sec:proof_maintheorem}
This section is devoted to the proof of Theorem \ref{thm:maintheorem}.
First, note that the factor $W(\mathcal{O}_\rho)$ in the statement is $\rho(G)$-invariant, so the product $W(\mathcal{O}_\rho)\langle G_\rho \rangle$ of two subgroups of $W \rtimes G$ is indeed the semidirect product $W(\mathcal{O}_\rho) \rtimes \langle G_\rho \rangle$.
This implies that, since $W(\mathcal{O}_\rho)$ is generated by involutions, the claim (2) follows immediately from (1).
So we prove (1) below.

For the ``only if'' part, we assume that $wg \in (W \rtimes G)_{\mathrm{ACI}}$ and prove that $w \in W(\mathcal{O}_\rho)$ and $g \in G_\rho \cup \{1\}$.
Now by (\ref{eq:index_intersection}) and Corollary \ref{cor:indexandZ} (2), we have
\begin{equation}
\label{eq:proofofmaintheorem_hypothesis}
\left[H:Z_H(w'g')\right]<\infty \textrm{ for any } H \leq W \rtimes G \textrm{ and } w'g' \in \langle wg \rangle_{\lhd W \rtimes G}.
\end{equation}
Note that $\rho_g(w)=w^{-1}$ and $g^2=1$ since $1=(wg)^2=w\rho_g(w) \cdot g^2$.
We divide the proof into the following five steps.\medskip\\
\textbf{Step 1: $\rho_g$ maps each $W_I \in \mathcal{C}_W^{\mathrm{inf}}$ onto itself.}\medskip

Assume contrary that $\rho_g$ maps $W_I$ onto an irreducible component other than $W_I$.
Let $\pi:W \twoheadrightarrow W_I$ be the projection.
Take $s \in I$ and put $a=swgs(wg)^{-1} \in \langle wg \rangle_{\lhd W \rtimes G}$.
Then we have $a=sw\rho_g(s)w^{-1} \in W$, so $Z_{W_I}(a)=Z_{W_I}(\pi(a))$.
Thus (\ref{eq:proofofmaintheorem_hypothesis}) implies that $\pi(a) \in W_I$ is almost central in $W_I$.
However, the first assumption yields that $\pi(\rho_g(s))=1$, so $\pi(a)=s\pi(w)1\pi(w)^{-1}=s$, which is not almost central in $W_I$ by Proposition \ref{prop:charoffinitepart}.
This is a contradiction.\medskip\\
\textbf{Step 2: $\rho_g$ is identity on every $W_I \in \mathcal{C}_W^{\mathrm{inf}}$.}\medskip

Assume that the claim fails for $W_I$.
Note that $\rho_g(W_I)=W_I$ by Step 1.
Let $\pi:W \twoheadrightarrow W_I$ be the projection.
Then we may assume without loss of generality that $\ell(\pi(w)) \leq \ell(\pi(uw\rho_g(u)^{-1}))$ for all $u \in W_I$; if this inequality fails, replace $wg$ with another involution $u(wg)u^{-1}=uw\rho_g(u)^{-1} \cdot g$ in $\langle wg \rangle_{\lhd W \rtimes G}$, which is also almost central in $W \rtimes G$ by (\ref{eq:proofofmaintheorem_hypothesis}), and use the induction on $\ell(\pi(w))$.

Put $\tau=\rho_g|_I \in \mathrm{Aut}\,\Gamma_I$, which is assumed to be non-identity.
Now if $\pi(w)=1$, then we have $Z_{W_I}(wg)=Z_{W_I}(g)={W_I}^\tau$ and so $\left[W_I:{W_I}^\tau\right]<\infty$ by (\ref{eq:proofofmaintheorem_hypothesis}), contradicting Theorem \ref{thm:indexoffixedpointsubgroup}.
Thus $\pi(w) \neq 1$.

We show that $\pi(w)$ is an involution in ${W_I}^\tau$.
Let $s_1 \cdots s_n$ (where $n \geq 1$ and $s_i \in I$) be an arbitrary reduced expression of $\pi(w) \in W_I$.
Then, since $\rho_g(w)=w^{-1}$ and $\rho_g(W_I)=W_I$, we have
\[
\pi(w)=\pi(\rho_g(w)^{-1})=\rho_g(\pi(w)^{-1})=\tau(\pi(w)^{-1})=\tau(s_n) \cdots \tau(s_1),
\]
so $\ell(\pi(w)\tau(s_1))<\ell(\pi(w))$, therefore Exchange Condition shows that $\pi(w)=s_1 \cdots \widehat{s_i} \cdots s_n\tau(s_1)$ for an index $i$.
Now if $i \geq 2$, then $\pi(s_1w\tau(s_1)^{-1})=s_2 \cdots \widehat{s_i} \cdots s_n$ is shorter than $\pi(w)$, contradicting the minimality of $\ell(\pi(w))$.
Thus we have $i=1$ and $\pi(w)=s_2 \cdots s_n\tau(s_1)$.
Since the original reduced expression $s_1 \cdots s_n$ is arbitrary, we can apply this argument to the new expression of $\pi(w)$.
Iterating, we have
\[
\pi(w)=s_3 \cdots s_n\tau(s_1)\tau(s_2)= \cdots =s_n\tau(s_1) \cdots \tau(s_{n-1})=\tau(s_1) \cdots \tau(s_n).
\]
Since $\pi(w)=s_1 \cdots s_n=\tau(s_n) \cdots \tau(s_1)$, the claim of this paragraph follows.

Now if $m_{s,\tau(s)}=\infty$ for some $s \in I$, then since $\tau^2=\mathrm{id}_I$, Theorem \ref{thm:Wtauisinfinite} (2) (applied to $\pi(w)$) gives us an element $u \in W_I$ of infinite order such that
\[
\pi(u^kwgu^{-k}(wg)^{-1})=u^k\pi(w)\tau(u)^{-k}\pi(w)^{-1} \neq 1 \textrm{ for all } k \neq 0.
\]
This means that $Z_{\langle u \rangle}(wg)=1$, so $\left[\!\right.\langle u \rangle:Z_{\langle u \rangle}(wg)\left.\!\right]=\infty$, contradicting (\ref{eq:proofofmaintheorem_hypothesis}).
On the other hand, if $m_{s,\tau(s)}<\infty$ for all $s \in I$, then the Coxeter group ${W_I}^\tau$ is infinite and irreducible by Theorem \ref{thm:Wtauisinfinite} (1).
Now we have
\[
Z_{{W_I}^\tau}(\pi(w))=Z_{{W_I}^\tau}(w)=Z_{{W_I}^\tau}(wg),
\]
which has finite index in ${W_I}^\tau$ by (\ref{eq:proofofmaintheorem_hypothesis}).
Thus the non-identity involution $\pi(w) \in {W_I}^\tau$ is almost central in ${W_I}^\tau$, contradicting Proposition \ref{prop:charoffinitepart}.
Hence Step 2 is concluded.\medskip\\
\textbf{Step 3: $w \in W_{\mathrm{fin}}$.}\medskip

We show that $\pi(w)=1$ for any $W_I \in \mathcal{C}_W^{\mathrm{inf}}$ with projection $\pi:W \twoheadrightarrow W_I$.
Since $\rho_g$ is identity on $W_I$ (Step 2), we have $Z_{W_I}(wg)=Z_{W_I}(w)=Z_{W_I}(\pi(w))$ and so (by (\ref{eq:proofofmaintheorem_hypothesis})) $\pi(w)$ is almost central in $W_I$.
Now since $1=\pi(w\rho_g(w))=\pi(w)\rho_g(\pi(w))=\pi(w)^2$, the claim follows from Proposition \ref{prop:charoffinitepart}.\medskip\\
\textbf{Step 4: $w \in W(\mathcal{O}_{\rho})$.}\medskip

Assume the contrary.
Then there exist a $\rho^\dagger(G)$-orbit $\mathcal{O} \subseteq \mathcal{C}_W^{\mathrm{fin}}$ with infinite cardinality and $W_I \in \mathcal{O}$ (with projection $\pi_I:W \twoheadrightarrow W_I$) such that $\pi_I(w) \neq 1$.
Fix the $\mathcal{O}$, and let $\mathcal{O}_0$ be the set of all such $W_I \in \mathcal{O}$, so $|\mathcal{O}_0|<\infty$.

We show that $\rho^\dagger_h(\mathcal{O}_0)=\mathcal{O}_0$, or equivalently $\rho^\dagger_h(\mathcal{O}_0) \subseteq \rho^\dagger_h(\mathcal{O}_0)$, for any $h \in Z_G(wg)$.
Note that $\rho_h(w)=w$ since $wgh=hwg=\rho_h(w)hg$.
Now if $W_I \in \mathcal{O}_0$ and $\rho^\dagger_h(W_I)=W_J \not\in \mathcal{O}_0$, then $\pi_J(\rho_h(w))=\rho_h(\pi_I(w)) \neq 1$ and $\pi_J(w)=1$, contradicting $\rho_h(w)=w$.
Thus $\rho^\dagger_h(\mathcal{O}_0) \subseteq \mathcal{O}_0$ as desired.

Since $\mathcal{O}$ is an infinite $\rho^\dagger(G)$-orbit, we can choose infinitely many finite subsets $\mathcal{O}_1,\mathcal{O}_2,\dots$ of $\mathcal{O}$, irreducible components $W_{I_0},W_{I_1},W_{I_2},\dots \in \mathcal{O}$ and elements $g_0,g_1,g_2,\dots \in G$ inductively, where we start with arbitrary $W_{I_0} \in \mathcal{O}_0$ and $g_0=1$, subject to the conditions $W_{I_k} \not\in \bigcup_{i=0}^{k-1}\mathcal{O}_i$, $\rho^\dagger_{g_k}(W_{I_0})=W_{I_k}$ and $\mathcal{O}_k=\rho^\dagger_{g_k}(\mathcal{O}_0) \ni W_{I_k}$ for all $k \geq 1$.
Now if $i<j$ and $h \in Z_G(wg)$, then the previous paragraph implies that $\rho^\dagger_{g_ih}(W_{I_0})=\rho^\dagger_{g_i}\rho^\dagger_h(W_{I_0}) \in \rho^\dagger_{g_i}(\mathcal{O}_0)=\mathcal{O}_i$, while $\rho^\dagger_{g_j}(W_{I_0})=W_{I_j} \not\in \mathcal{O}_i$, so we have $g_j \neq g_ih$.
Thus all the $g_i$ belong to distinct cosets in $G/Z_G(wg)$, while $\left[G:Z_G(wg)\right]<\infty$ by (\ref{eq:proofofmaintheorem_hypothesis}).
This contradiction yields the claim.\medskip\\
\textbf{Step 5: $g \in G_\rho \cup \{1\}$.}\medskip

Note that $g^2=1$.
Since $hwg=\rho_h(w) \cdot hg$ for $h \in G$, we have $Z_G(wg) \subseteq Z_G(g)$ and so $g$ is almost central in $G$ by (\ref{eq:proofofmaintheorem_hypothesis}).
From now, we check (\ref{eq:condition_maintheorem}).

By Step 4, the union $\mathcal{O}$ of a finite number of some $\rho^\dagger(G)$-orbits with finite cardinalities satisfies that $w \in W(\mathcal{O})$.
Since $|\mathcal{O}|<\infty$, it suffices to show that $\rho_g$ is identity on all $W_I \in \mathcal{O}'=\mathcal{C}_W \smallsetminus \mathcal{O}$ except a finite number of finite irreducible components.
Now $\mathcal{O}'$ is $\rho^\dagger_g$-invariant as well as its complement $\mathcal{O}$, while $W(\mathcal{O}') \subseteq Z_W(w)$, so $Z_{W(\mathcal{O}')}(wg)=Z_{W(\mathcal{O}')}(g)$ is the fixed-point subgroup $W(\mathcal{O}')^{\tau}$ (where $\tau=\rho_g|_{W(\mathcal{O}')}$).
Since $Z_{W(\mathcal{O}')}(wg)$ has finite index in $W(\mathcal{O}')$ by (\ref{eq:proofofmaintheorem_hypothesis}), the claim follows from Theorem \ref{thm:indexoffixedpointsubgroup}.\medskip

Hence the ``only if'' part has been proved.
From now, we prove the other part; so we assume that $w \in W(\mathcal{O}_\rho)$, $g \in G_\rho \cup \{1\}$ and $wg$ is an involution, and prove that $wg$ is almost central in $W \rtimes G$.
By the choice of $w$, there are a finite number of finite $\rho^\dagger(G)$-orbits in $\mathcal{C}_W^{\mathrm{fin}}$ such that their union $\mathcal{O}$ satisfies that $w \in W(\mathcal{O})$.
Now note that
\[
Z_{W \rtimes G}(wg) \supseteq Z_{W \rtimes G}(w) \cap Z_{W \rtimes G}(g) \supseteq (Z_W(w)Z_G(w)) \cap (Z_W(g)Z_G(g)),
\]
so it suffices to show that both $Z_W(w)Z_G(w)$ and $Z_W(g)Z_G(g)$ have finite index in $W \rtimes G$ (see (\ref{eq:index_diamond})).
Moreover, Lemma \ref{lem:indexinsemidirect} reduces the claim to the following four claims:\medskip\\
\textbf{Step 6: $Z_W(w)$ has finite index in $W$.}\medskip

This follows since $w$ lies in the finite direct factor $W(\mathcal{O})$ of $W$.\medskip\\
\textbf{Step 7: $Z_G(w)$ has finite index in $G$.}\medskip

Since $W(\mathcal{O})$ is finite and $\rho(G)$-invariant, the action gives rise to a homomorphism $\rho'$ from $G$ to the finite group $\mathrm{Aut}\,W(\mathcal{O})$.
Now $\ker \rho'$ is contained in $Z_G(w)$ (since $w \in W(\mathcal{O})$) and has finite index in $G$, proving the claim.\medskip\\
\textbf{Step 8: $Z_W(g)$ has finite index in $W$.}\medskip

This is trivial if $g=1$.
If $g \in G_\rho$, then the property (\ref{eq:condition_maintheorem}) and Theorem \ref{thm:indexoffixedpointsubgroup} imply that the fixed-point subgroup $W^{\rho_g}=Z_W(g)$ by $\rho_g$ has finite index in $W$, as desired.\medskip\\
\textbf{Step 9: $Z_G(g)$ has finite index in $G$.}\medskip

This is obvious from the choice of $g$.\medskip

Hence the proof of Theorem \ref{thm:maintheorem} is concluded.

\noindent
\textbf{Koji Nuida}\\
Graduate School of Mathematical Sciences, University of Tokyo\\
\indent 3-8-1 Komaba, Meguro-ku, Tokyo, 153-8914, Japan\\
E-mail: nuida@ms.u-tokyo.ac.jp

\end{document}